\documentclass[11pt,draft]{article}
\usepackage{amsmath,amscd}
\usepackage{amssymb,latexsym,amsthm}
\usepackage{color}
\usepackage[spanish,english]{babel}
\usepackage{amssymb}
\usepackage{color}
\usepackage{amsmath,amsthm,amscd}
\usepackage[latin1]{inputenc}
\usepackage{lscape}
\usepackage{amsfonts}

\numberwithin{equation}{section}
\newtheorem{theorem}{Theorem}[section]

\newtheorem{proposition}[theorem]{Proposition}
\newtheorem{lemma}[theorem]{Lemma}

\newtheorem{corollary}[theorem]{Corollary}

\theoremstyle{definition}

\newtheorem{definition}[theorem]{Definition}
\newtheorem{examples}[theorem]{Examples}

\newtheorem{remark}[theorem]{Remark}

\newcommand{\cB}{\mbox{${\cal B}$}}

\newcommand{\cS}{\mbox{${\cal S}$}}

\newcommand{\cU}{\mbox{${\cal U}$}}

\newcommand{\cW}{\mbox{${\cal W}$}}

\title{\textbf{$\sigma$-PBW Extensions of Skew Armendariz Rings}}
\author{Armando Reyes\footnote{Departamento de  Matem\'aticas. e-mail: mareyesv@unal.edu.co}\\ Universidad Nacional de Colombia \\ H\'ector Su\'arez\footnote{Escuela de Matem\'aticas y Estad\'istica. e-mail: hjsuarezs@unal.edu.co}\\ Universidad Pedag\'ogica y Tecnol\'ogica de  Colombia}
\date{}
\begin{document}
\maketitle
\begin{abstract}

\noindent The aim of this paper is to investigate a general notion of $\sigma$-PBW extensions over Armendariz rings. As an application,  the properties of being Baer, quasi-Baer, p.p. and p.q.-Baer are established for these extensions. We generalize several results in the literature for Ore extensions of injective type and skew PBW extensions.

\bigskip


\bigskip




\end{abstract}

\section{Introduction}
A commutative ring $B$ is called {\em Armendariz} (the term was introduced by Rege and Chhawchharia in \cite{RegeChhawchharia1997}) if for polynomials $f(x)=a_0+a_1x+\dotsb + a_nx^n$, $g(x)=b_0+b_1x+\dotsb + b_mx^m$ of $B[x]$ which satisfy $f(x)g(x)=0$, then $a_ib_j=0$, for every $i,j$. As we can appreciate, the interest of this notion lies in its natural and its useful role in understanding the relation between the annihilators of the ring $B$ and the annihilators of the polynomial ring $B[x]$. For instance, in   \cite{Armendariz1974}, Lemma 1, Armendariz showed that a {\em reduced} ring (i.e., a ring without nonzero nilpotent elements) always satisfies this condition (reduced rings are {\em Abelian} that is, every idempotent is central, and also semiprime, i.e., its prime radical is trivial). For non-commutative rings, more exactly the well-known Ore extensions introduced by Ore in \cite{Ore1933}, the notion of Armendariz has been also studied. Commutative and non-commutative treatments have been investigated in several papers, see \cite{Armendariz1974}, \cite{RegeChhawchharia1997}, \cite{AndersonCamillo1998},  \cite{KimLee2000}, \cite{Huhetal2002}, \cite{HongKimKwak2003}, \cite{LeeWong2003}, \cite{Matczuk2004}, and others.\\

The non-commutative rings of interest for us in this article are the $\sigma$-PBW extensions (also known as skew Poincar\'e-Birkhoff-Witt extensions) introduced in \cite{LezamaGallego}. These structures are more general than iterated Ore extensions of injective type defined by Ore in \cite{Ore1933}, universal enveloping algebras of finite dimensional Lie algebras, PBW extensions introduced by Bell and Goodearl in \cite{BellGoodearl1988}, almost normalizing extensions defined by McConnell and Robson in \cite{McConnellRobson2001}, solvable polynomial rings introduced by Kandri-Rody and Weispfenning in \cite{KandryWeispfenninig1990}, and generalized by Kredel in \cite{Kredel1993},  diffusion algebras introduced by Isaev, Pyatov, and Rittenberg  in \cite{IsaevPyatovRittenberg2001}, and other algebras. The importance of $\sigma$-PBW extensions is that we do not assume that the coefficients commute with the variables, and we take coefficients not necessarily in fields (see Definition \ref {gpbwextension}). In fact, the $\sigma$-PBW extensions contain well-known groups of algebras such as some types of $G$-algebras in the sense of Levandovskyy \cite{Levandovskyy2005}, Auslander-Gorenstein rings, some Calabi-Yau and skew Calabi-Yau algebras, some Artin-Schelter regular algebras, some Koszul algebras, quantum polynomials, some quantum universal enveloping algebras, and others (see \cite{Reyes2013PhD} or
\cite{Lezama-Reyes1} for a detailed list of examples). With respect to Clifford and Grassman algebras, skew PBW extensions are also important. More precisely, since any Clifford algebra is a quotient of a sol\-va\-ble polynomial ring by a two-sided ideal (\cite{KandryWeispfenninig1990}, p. 24), and solvable polynomial rings are strictly contained in $\sigma$-PBW extensions, then Cliford and Grassman algebras are quotients of $\sigma$-PBW extensions by two-sided ideals (in \cite{Reyes2013PhD} and \cite{LezamaAcostaReyes2015} it was  presented a characterization of ideals in $\sigma$-PBW extensions). For more details about the relation between $\sigma$-PBW extensions and another algebras with PBW bases, see \cite{Reyes2013PhD} or
\cite{Lezama-Reyes1}.\\

Now, since ring, module and homological properties of $\sigma$-skew PBW extensions have been studied by the authors and others (see \cite{LezamaGallego}, \cite{Lezamaetal2013}, \cite{Reyes2013PhD}, \cite{Reyes2013}, \cite{Reyes2014}, \cite{LezamaAcostaReyes2015},  \cite{Reyes2015}, \cite{SuarezLezamaReyes2015}, \cite{NinoReyes2016}, \cite{ReyesSuarez2016a}, \cite{ReyesSuarez2016C},  \cite{ReyesSuarezUPTC2016}, and others), the aim of this paper is to establish a general notion of {\em skew} Armendariz ring for $\sigma$-PBW extensions which generalize the case of Ore extensions and previous papers, and formulate new results for several non-commutative algebras which can not be expressed as iterated Ore extensions. The theory developed here generalizes the treatments presented in \cite{NasrMoussavi2008} for Ore extensions of injective type, and the results established in \cite{Reyes2015}, \cite{NinoReyes2016}, and \cite{ReyesSuarez2016C} for skew PBW extensions. As an application of our treatment, we characterize the properties of being Baer, quasi-Baer, p.p. and p.q.-Baer for these extensions generalizing several results in the literature for Ore extensions of injective type and skew PBW extensions (\cite{HongKimKwak2000}, \cite{KimLee2000},  \cite{Hirano2002}, \cite{HongKimKwak2003}, \cite{MoussaviHashemi2005},  \cite{MohamedMoufadal2006}, \cite{NasrMoussavi2007}, \cite{NasrMoussavi2008}, \cite{HongKimLee2009}, \cite{Reyes2015}, \cite{NinoReyes2016}, and \cite{ReyesSuarez2016C}).\\

Next, we describe the structure of this article. In Section \ref{definitionexamplesspbw} we establish some useful results about $\sigma$-PBW extensions for the rest of the paper. In Section \ref{reviewwww} we present a review of the proposals about a notion of Armendariz ring for these extensions. Then, in Section \ref{skewArmendarizsection} we introduce two notions of Armendariz: skew-Armendariz (Definition \ref{perritogatito}) and a more general notion, the weak skew-Armendariz (Definition \ref{Definition2.42008}). These definitions are motivated by the treatment de\-ve\-lo\-ped in \cite{NasrMoussavi2008} for Ore extensions, and generalize the theory  for classical polynomial rings, Ore extensions of injective type, and skew PBW extensions (c.f. \cite{Armendariz1974}, \cite{RegeChhawchharia1997}, \cite{AndersonCamillo1998}, \cite{KimLee2000}, \cite{Huhetal2002}, \cite{HongKimKwak2003}, \cite{LeeWong2003}, \cite{Matczuk2004}, and \cite{NasrMoussavi2007}, \cite{Reyes2015}, \cite{NinoReyes2016}, \cite{ReyesSuarez2016C}, and others). In particular, we show that the families  of Armendariz rings defined in \cite{NinoReyes2016} and \cite{ReyesSuarez2016C} are strictly contained in the family of skew-Armendariz (or the family of weak skew-Armendariz) but the converse is false. In this section we also prove that if $R$ is a weak skew-Armendariz ring, then $R$ and $A$ are Abelian (Propositions \ref{Theorem3.32008} and \ref{Corollary3.42008}, respectively), and we characterize the property of skew-Armendariz over the classical right quotient rings of $R$ (Theorem \ref{Theorem4.72008}). In Section \ref{QuasiBaerrings} we generalize the results presented by the authors in \cite{Reyes2015}, \cite{NinoReyes2016}, and \cite{ReyesSuarez2016C}, about the characterizations of being Baer, quasi-Baer, p.p. and the p.q.-Baer for $\sigma$-PBW extensions. We adapt the ideas presented in \cite{NasrMoussavi2008} and use the notion of $(\Sigma,\Delta)$-ideal (this notion was used in \cite{Reyes2014} and \cite{LezamaAcostaReyes2015} for the study of the uniform dimension and the prime ideals of $\sigma$-PBW extensions, respectively), and the notion of $(\Sigma, \Delta)$-quasi-baer rings. \\

Throughout the paper, the word ring means a ring with unity not necessarily commutative.

\section{Definitions and elementary properties}\label{definitionexamplesspbw}
We start recalling the definition and some preliminaries about our object of study.
\begin{definition}[\cite{LezamaGallego}, Definition 1]\label{gpbwextension}
Let $R$ and $A$ be rings. We say that $A$ {\rm is a $\sigma$-PBW extension of} $R$ (also called  \textit{skew PBW
extension of} $R$), if the following conditions hold:
\begin{enumerate}
\item[\rm (i)]$R\subseteq A$;
\item[\rm (ii)]there exist elements $x_1,\dots ,x_n\in A$ such that $A$ is a left free $R$-module with basis  ${\rm Mon}(A):= \{x^{\alpha}=x_1^{\alpha_1}\cdots
x_n^{\alpha_n}\mid \alpha=(\alpha_1,\dots ,\alpha_n)\in
\mathbb{N}^n\}$, and $x_1^{0}\dotsb x_n^{0}:=1\in {\rm Mon}(A)$.
\item[\rm (iii)]For each $1\leq i\leq n$ and any $r\in R\ \backslash\ \{0\}$, there exists an element $c_{i,r}\in R\ \backslash\ \{0\}$ such that $x_ir-c_{i,r}x_i\in R$.
\item[\rm (iv)]For any elements $1\leq i,j\leq n$, there exists $c_{i,j}\in R\ \backslash\ \{0\}$ such that $x_jx_i-c_{i,j}x_ix_j\in R+Rx_1+\cdots +Rx_n$ (i.e., there exist elements $r_0^{(i,j)}, r_1^{(i,j)}, \dotsc, r_n^{(i,j)} \in R$ with $x_jx_i - c_{i,j}x_ix_j = r_0^{(i,j)} + \sum_{k=1}^{n} r_k^{(i,j)}x_k$).
\end{enumerate}
If all these conditions are satisfied, we will write $A:=\sigma(R)\langle
x_1,\dots,x_n\rangle$.
\end{definition}
\begin{proposition}[\cite{LezamaGallego}, Proposition
3]\label{sigmadefinition}
Let $A$ be a skew PBW extension of $R$. For each $1\leq i\leq
n$, there exist an injective endomorphism $\sigma_i:R\rightarrow
R$ and a $\sigma_i$-derivation $\delta_i:R\rightarrow R$, such that $x_ir=\sigma_i(r)x_i+\delta_i(r)$, for  each $r\in R$. We write $\Sigma:=\{\sigma_1,\dotsc, \sigma_n\}$, and $\Delta:=\{\delta_1,\dotsc, \delta_n\}$.
\end{proposition}
\begin{remark}\label{notesondefsigampbw}
With respect to the Definition \ref{gpbwextension} and the Proposition \ref{sigmadefinition}, we have the following remarks:
\begin{itemize}
\item Since ${\rm Mon}(A)$ is a left $R$-basis of $A$, the elements $c_{i,r}$ and $c_{i,j}$ in Definition \ref{gpbwextension} are unique.
\item In Definition \ref{gpbwextension} (iv), $c_{i,i}=1$. This follows from the equality $x_i^2-c_{i,i}x_i^2=s_0+s_1x_1+\cdots+s_nx_n$, with $s_i\in R$, which implies $1-c_{i,i}=0=s_i$.
\item If $i<j$ and $d_i', b_j'\in R$, then $d_i'x_ib_j'x_j = d_i'[\sigma_i(b_j')x_i + \delta_i(b_j')]x_j = d_i'\sigma_i(b_j')x_ix_j + d_i'\delta_i(b_j')x_j$. Since $x_jx_i = c_{i,j}x_ix_j + r^{(i,j)} + \sum_{k=1}^{n} r^{(i,j)}_kx_k$, then $d_j'x_jb_i'x_i = d_j'[\sigma_j(b_i')x_j + \delta_j(b_i')]x_i = d_j'\sigma_j(b_i')x_jx_i$\ $+\ d_j'\delta_j(b_i')x_i =$ $  d_j'\sigma_j(b_i')(c_{i,j}x_ix_j + r^{(i,j)} + \sum_{k=1}^{n} r^{(i,j)}_kx_k) + d_j'\delta_j(b_i')x_i$. In this way,
\begin{align*}
d_i'x_ib_j'x_j + d_j'x_jb_i'x_i = &\ d_i'\sigma_i(b_j')x_ix_j + d_i'\delta_i(b_j')x_j\\
+ &\ d_j'\sigma_j(b_i')\biggl(c_{i,j}x_ix_j + r^{(i,j)} + \sum_{k=1}^{n} r^{(i,j)}_kx_k\biggr) + d_j'\delta_j(b_i')x_i\\
= &\ [d_i'\sigma_i(b_j') + d_j'\sigma_j(b_i')c_{i,j}]x_ix_j + d_j'\delta_j(b_i')x_i \\
+ &\ d_i'\delta_i(b_j')x_j
+ d_j'\sigma_j(b_i')r^{(i,j)} + d_j'\sigma_j(b_i')\sum_{k=1}^{n} r^{(i,j)}_kx_k\\
= &\ [d_i'\sigma_i(b_j') + d_j'\sigma_j(b_i')c_{i,j}]x_ix_j + [d_j'\delta_j(b_i') + d_j'\sigma_j(b_i')r_i^{(i,j)}]x_i \\
+ &\ [d_i'\delta_i(b_j') + d_j'\sigma_j(b_i')r_j^{(i,j)}]x_j
+ d_j'\sigma_j(b_i')r^{(i,j)} + d_j'\sigma_j(b_i')\sum_{k=1,\ k\neq i, j}^{n} r^{(i,j)}_kx_k.
\end{align*}
\end{itemize}
\end{remark}
\begin{definition}(\cite{LezamaGallego}, Definition 4)\label{sigmapbwderivationtype}
Consider $A$ a skew PBW extension of $R$.
\begin{enumerate}
\item[\rm (i)] $A$ is called \textit{quasi-commutative} if the conditions
{\rm(}iii{\rm)} and {\rm(}iv{\rm)} in Definition
\ref{gpbwextension} are replaced by the following: (iii') for each $1\leq i\leq n$ and all $r\in R\ \backslash\ \{0\}$, there exists $c_{i,r}\in R\ \backslash\ \{0\}$ such that $x_ir=c_{i,r}x_i$; (iv') for any $1\leq i,j\leq n$, there exists $c_{i,j}\in R\ \backslash\ \{0\}$ such that $x_jx_i=c_{i,j}x_ix_j$.
\item[\rm (ii)] $A$ is called \textit{bijective} if $\sigma_i$ is bijective for each $1\leq i\leq n$, and $c_{i,j}$ is invertible, for any $1\leq
i<j\leq n$.
\end{enumerate}
\end{definition}
\begin{examples}\label{mentioned}
If $R[x_1;\sigma_1,\delta_1]\dotsb [x_n;\sigma_n,\delta_n]$ is an iterated Ore extension where
\begin{itemize}
\item $\sigma_i$ is injective, for $1\le i\le n$;
\item $\sigma_i(r)$, $\delta_i(r)\in R$, for every $r\in R$ and $1\le i\le n$;
\item $\sigma_j(x_i)=cx_i+d$, for $i < j$, and $c, d\in R$, where $c$ has a left inverse;
\item $\delta_j(x_i)\in R + Rx_1 + \dotsb + Rx_n$, for $i < j$,
\end{itemize}
then $R[x_1;\sigma_1,\delta_1]\dotsb [x_n;\sigma_n, \delta_n] \cong \sigma(R)\langle x_1,\dotsc, x_n\rangle$ (\cite{Lezama-Reyes1}, p. 1212). In particular, note that skew PBW extensions of endomorphism type are more general than iterated Ore extensions $R[x_1;\sigma_1]\dotsb [x_n;\sigma_n]$. On the other hand, skew PBW extensions are more general than Ore extensions of injective type (diffusion algebras, universal enveloping algebras of finite Lie algebras, and others, are examples of skew PBW extensions which can not be expressed as iterated Ore extensions, see \cite{Lezama-Reyes1} for more details). Skew PBW extensions contains various well-known groups of algebras such as PBW extensions \cite{BellGoodearl1988}, the almost normalizing extensions  \cite{McConnellRobson2001}, solvable polynomial rings \cite{KandryWeispfenninig1990}, and  \cite{Kredel1993}, diffusion algebras \cite{IsaevPyatovRittenberg2001}, some types of Auslander-Gorenstein rings, some skew Calabi-Yau algebras, some Artin-Schelter regular algebras, some Koszul algebras, quantum polynomials, some quantum universal enveloping algebras,  etc. In comparison with $G$-algebras \cite{Levandovskyy2005}, $\sigma$-PBW extensions do not assume that the ring of coefficients is a field neither that the coefficients commute with the variables, so that skew PBW extensions are not included in these algebras. Indeed, the $G$-algebras with  $d_{i,j}$ linear (recall that for these algebras $x_jx_i = c_{i,j}x_ix_j+ d_{i,j},\ 1\le i < j \le n$), are particular examples of $\sigma$-PBW extensions.  A detailed list of examples of skew PBW extensions and its relations with another algebras with PBW bases is presented in \cite{Reyes2013PhD},  \cite{Lezama-Reyes1}, and \cite{SuarezLezamaReyes2015}.
\end{examples}
\begin{definition}[\cite{LezamaGallego}, Definition 6]\label{definitioncoefficients}
Let $A$ be a skew PBW extension of $R$ with endomorphisms
$\sigma_i$, $1\leq i\leq n$, as in Proposition
\ref{sigmadefinition}.
\begin{enumerate}
\item[\rm (i)]For $\alpha=(\alpha_1,\dots,\alpha_n)\in \mathbb{N}^n$,
$\sigma^{\alpha}:=\sigma_1^{\alpha_1}\cdots \sigma_n^{\alpha_n}$,
$|\alpha|:=\alpha_1+\cdots+\alpha_n$. If
$\beta=(\beta_1,\dots,\beta_n)\in \mathbb{N}^n$; then
$\alpha+\beta:=(\alpha_1+\beta_1,\dots,\alpha_n+\beta_n)$.
\item[\rm (ii)]For $X=x^{\alpha}\in {\rm Mon}(A)$,
$\exp(X):=\alpha$, $\deg(X):=|\alpha|$, and $X_0:=1$. The symbol $\succeq$ will denote a total order defined on ${\rm Mon}(A)$ (a total order on $\mathbb{N}^n$). For an
 element $x^{\alpha}\in {\rm Mon}(A)$, ${\rm exp}(x^{\alpha}):=\alpha\in \mathbb{N}^n$.  If
$x^{\alpha}\succeq x^{\beta}$ but $x^{\alpha}\neq x^{\beta}$, we
write $x^{\alpha}\succ x^{\beta}$. Every element $f\in A$ can be expressed uniquely as $f=a_0 + a_1X_1+\dotsb +a_mX_m$, with $a_i\in R$, and $X_m\succ \dotsb \succ X_1$. With this notation, we define ${\rm
lm}(f):=X_m$, the \textit{leading monomial} of $f$; ${\rm
lc}(f):=a_m$, the \textit{leading coefficient} of $f$; ${\rm
lt}(f):=a_mX_m$, the \textit{leading term} of $f$; ${\rm exp}(f):={\rm exp}(X_m)$, the \textit{order} of $f$; and
 $E(f):=\{{\rm exp}(X_i)\mid 1\le i\le t\}$. Note that $\deg(f):={\rm max}\{\deg(X_i)\}_{i=1}^t$. Finally, if $f=0$, then
${\rm lm}(0):=0$, ${\rm lc}(0):=0$, ${\rm lt}(0):=0$. We also
consider $X\succ 0$ for any $X\in {\rm Mon}(A)$. For a detailed description of monomial orders in skew PBW  extensions, see \cite{LezamaGallego}, Section 3.
\end{enumerate}
\end{definition}
\begin{proposition}[\cite{LezamaGallego}, Theorem 7]\label{coefficientes}
If $A$ is a polynomial ring with coefficients in $R$ and the set of variables $\{x_1,\dots,x_n\}$, then $A$ is a skew PBW  extension of $R$ if and only if the following conditions are satisfied:
\begin{enumerate}
\item[\rm (i)]for each $x^{\alpha}\in {\rm Mon}(A)$ and every $0\neq r\in R$, there exist unique elements $r_{\alpha}:=\sigma^{\alpha}(r)\in R\ \backslash\ \{0\}$, $p_{\alpha ,r}\in A$, such that $x^{\alpha}r=r_{\alpha}x^{\alpha}+p_{\alpha, r}$,  where $p_{\alpha ,r}=0$ or $\deg(p_{\alpha ,r})<|\alpha|$, if
$p_{\alpha , r}\neq 0$. If $r$ is left invertible,  so is $r_\alpha$.
\item[\rm (ii)] For each $x^{\alpha},x^{\beta}\in {\rm Mon}(A)$ there exist unique elements $c_{\alpha,\beta}\in R$ and $p_{\alpha,\beta}\in A$ such that $x^{\alpha}x^{\beta}=c_{\alpha,\beta}x^{\alpha+\beta}+p_{\alpha,\beta}$, where $c_{\alpha,\beta}$ is left invertible, $p_{\alpha,\beta}=0$
or $\deg(p_{\alpha,\beta})<|\alpha+\beta|$, if
$p_{\alpha,\beta}\neq 0$.
\end{enumerate}
\end{proposition}
\begin{remark}(\cite{Reyes2015}, Remark 2.10)\label{juradpr}
If $X_i:=x_1^{\alpha_{i1}}\dotsb x_n^{\alpha_{in}}$ and $Y_j:=x_1^{\beta_{j1}}\dotsb x_n^{\beta_{jn}}$, then
\begin{align*}
a_iX_ib_jY_j = &\ a_i\sigma^{\alpha_i}(b_j)x^{\alpha_i}x^{\beta_j} + a_ip_{\alpha_{i1}, \sigma_{i2}^{\alpha_{i2}}(\dotsb (\sigma_{in}^{\alpha_{in}}(b)))} x_2^{\alpha_{i2}}\dotsb x_n^{\alpha_{in}}x^{\beta_j} \\
+ &\ a_i x_1^{\alpha_{i1}}p_{\alpha_{i2}, \sigma_3^{\alpha_{i3}}(\dotsb (\sigma_{{in}}^{\alpha_{in}}(b)))} x_3^{\alpha_{i3}}\dotsb x_n^{\alpha_{in}}x^{\beta_j} \\
+ &\ a_i x_1^{\alpha_{i1}}x_2^{\alpha_{i2}}p_{\alpha_{i3}, \sigma_{i4}^{\alpha_{i4}} (\dotsb (\sigma_{in}^{\alpha_{in}}(b)))} x_4^{\alpha_{i4}}\dotsb x_n^{\alpha_{in}}x^{\beta_j}\\
+ &\ \dotsb + a_i x_1^{\alpha_{i1}}x_2^{\alpha_{i2}} \dotsb x_{i(n-2)}^{\alpha_{i(n-2)}}p_{\alpha_{i(n-1)}, \sigma_{in}^{\alpha_{in}}(b)}x_n^{\alpha_{in}}x^{\beta_j} \\
+ &\ a_i x_1^{\alpha_{i1}}\dotsb x_{i(n-1)}^{\alpha_{i(n-1)}}p_{\alpha_{in}, b}x^{\beta_j}.
\end{align*}
Using this expression, we can see that when we compute every summand of $a_iX_ib_jY_j$ we obtain products of the coefficient $a_i$ with several evaluations of $b_j$ in $\sigma$'s and $\delta$'s depending of the coordinates of $\alpha_i$.
\end{remark}
\section{A review of the Armendariz notions for $\sigma$-PBW extensions}\label{reviewwww}
For a ring $B$ with a ring endomorphism $\sigma:B\to B$, and a $\sigma$-derivation $\delta:B\to B$, Krempa in \cite{Krempa1996} considered the Ore extension $B[x;\sigma,\delta]$, and defined $\sigma$ as a  {\em rigid en\-do\-mor\-phism} if $b\sigma(b)=0$ implies $b=0$, for $b\in B$. Krempa called $B$ $\sigma$-{\em rigid} if there exists a rigid endomorphism $\sigma$ of $B$. Properties of being Baer, quasi-Baer, p.p., and p.q.-Baer over $\sigma$-rigid rings have been investigated (c.f. \cite{Krempa1996}, \cite{HongKimKwak2000}, and others). All these results were generalized by the first author to the class of $\sigma$-PBW extensions. There, the key fact was an adequate notion of {\em rigidness} for these extensions. Let us recall it.
\begin{definition}
(\cite{Reyes2015}, Definition 3.2)
If $B$ is a ring and $\Sigma$ a family of endomorphisms of $B$, then $\Sigma$ is called a {\em rigid endomorphisms family} if $r\sigma^{\alpha}(r)=0$ implies $r=0$, for every $r\in B$ and $\alpha\in \mathbb{N}^n$. A ring $B$ is called to be $\Sigma$-{\em rigid} if there exists a rigid endomorphisms family $\Sigma$ of $B$.
\end{definition}

Note that if $\Sigma$ is a rigid endomorphisms family, then every element $\sigma_i\in \Sigma$ is a monomorphism. In this way, we consider the family of injective endomorphisms $\Sigma$ and the family $\Delta$ of $\Sigma$-derivations of a skew PBW extension $A$ of a ring $R$ (see Proposition \ref{sigmadefinition}).  $\Sigma$-rigid rings are reduced rings: if $B$ is a $\Sigma$-rigid ring and $r^2=0$, for $r\in B$, then $0=r\sigma^{\alpha}(r^2)\sigma^{\alpha}(\sigma^{\alpha}(r))=r\sigma^{\alpha}(r)\sigma^{\alpha}(r)\sigma^{\alpha}(\sigma^{\alpha}(r))=r\sigma^{\alpha}(r)\sigma^{\alpha}(r\sigma^{\alpha}(r))$, i.e., $r\sigma^{\alpha}(r)=0$, and so $r=0$, that is, $B$ is reduced. By \cite{Reyes2015}, Corollary 3.4, if $A$ is a skew PBW  extension of a $\Sigma$-rigid ring $R$, the equality $ab=0$, for $a,b\in R$, implies $ax^{\alpha}bx^{\beta}=0$ in $A$, for any $\alpha, \beta\in \mathbb{N}^n$. Recall that if $A$ is a $\sigma$-PBW extension of a ring $R$, then $R$ is $\Sigma$-rigid if and only if $A$ is a reduced ring (\cite{Reyes2015}, Proposition 3.5).\\

With the purpose of generalizing the notion of $\sigma$-rigid ring and studying the properties of being Baer, quasi-Baer, p.p. and p.q.-Baer (the notion of Armendariz ring also was studied) over this more general structure, several notions of {\em skew} Armendariz rings have been established in the literature  (c.f. \cite{HongKimKwak2003}, \cite{MoussaviHashemi2005}, \cite{NasrMoussavi2007}). Precisely, in \cite{NinoReyes2016} and  \cite{ReyesSuarez2016C}, the authors generalized all these treatments for the case of skew PBW extensions by introducing the following definitions.
\begin{definition}(\cite{NinoReyes2016}, Definitions 3.4 and 3.5)\label{nino}
Let $A$ be a skew PBW extension of a ring $R$. We say that $R$ is an $(\Sigma, \Delta)$-{\em Armendariz ring}, if for polynomials $f = a_0 + a_1X_1 + \dotsb + a_mX_m$ and $g= b_0 + b_1Y_1 + \dotsb + b_tY_t$ in $A$, the equality $fg=0$ implies $a_iX_ib_jY_j=0$, for every $i, j$. We say that $R$ is an $(\Sigma, \Delta)$-{\em weak Armendariz ring}, if for linear polynomials $f = a_0 + a_1x_1 + \dotsb + a_nx_n$ and $g= b_0 + b_1x_1 + \dotsb + b_nx_n$ in $A$, the equality $fg=0$ implies $a_ix_ib_jx_j=0$, for every $i, j$.
\end{definition}
Note that every $\Sigma$-rigid ring is a $(\Sigma, \Delta)$-skew Armendariz ring (\cite{NinoReyes2016}, Proposition 3.6).

\begin{definition}(\cite{ReyesSuarez2016C}, Definitions 3.1 and 3.2)\label{defiskewAr}
Let $A$ be a skew PBW extension of a ring $R$. $R$ is called a {\em $\Sigma$-skew Armendariz} ring, if for elements $f=\sum_{i=0}^{m} a_iX_i$ and $g = \sum_{j=0}^{t} b_jY_j$ in $A$, the equality $fg=0$ implies $a_i\sigma^{\alpha_i}(b_j)=0$,  for all $0\le i\le m$ and $0\le j\le t$, where $\alpha_i={\rm exp}(X_i)$. $R$ is called a {\em weak $\Sigma$-skew Armendariz} ring, if for elements $f=\sum_{i=0}^{n} a_ix_i$ and $g = \sum_{j=0}^{n} b_jx_j$ in $A$ ($x_0:=1$), the equality $fg=0$ implies $a_i\sigma_i(b_j)=0$, for all $0\le i, j \le n$ ($\sigma_0:={\rm id}_R$).
\end{definition}
Note that every $\Sigma$-skew Armendariz ring is a weak $\Sigma$-skew Armendariz ring. If $A$ is a skew PBW extension of a ring $R$, and if $R$ is $\Sigma$-rigid, then $R$ is $\Sigma$-skew Armendariz (\cite{ReyesSuarez2016C}, Proposition 3.4). The converse of this proposition is false as the following remark shows.
\begin{remark}\label{conversws}
\begin{itemize}
\item Consider the ring  $ B=\biggl\{\begin{pmatrix} a & t\\ 0 & a \end{pmatrix}\mid a\in \mathbb{Z},\ t\in \mathbb{Q}\biggr\}$.
Then $B$ is a commutative ring, and if we consider the automorphism $\sigma$ of $R$ given by $\sigma\biggl(\begin{pmatrix} a & t \\ 0 & a \end{pmatrix}\biggr) = \begin{pmatrix} a & t/2 \\ 0 & a \end{pmatrix}$. In \cite{HongKimKwak2003}, Example 1, it was shown that $R$ is $\sigma$-skew Armendariz and is not a $\sigma$-rigid. Since $\Sigma$-rigid and $\Sigma$-skew Armendariz are generalizations of $\sigma$-rigid and $\sigma$-skew Armendariz, respectively, this example shows that there exist an example of a $\Sigma$-skew Armendariz ring which is not $\Sigma$-rigid.
\item Let $B=\mathbb{Z}_2[x]$ be the commutative polynomial ring over $\mathbb{Z}_2$, and $\sigma$ the endomorphism of $B=\mathbb{Z}_2[x]$ defined by $\sigma(f(x)) = f(0)$. Then $B=\mathbb{Z}_2[x]$ is $\sigma$-skew Armendariz and is not $\sigma$-rigid (\cite{HongKimKwak2003}, Example 5).
\end{itemize}
\end{remark}
From definitions above we can establish the following relations
\begin{equation*}
\Sigma{\rm -rigid}\ \subsetneqq\ (\Sigma, \Delta){\rm -Armendariz} \ \subsetneqq\ (\Sigma, \Delta){\rm -weak\ Armendariz}
\end{equation*}
\begin{equation*}
\Sigma{\rm-rigid}\ \subsetneqq\ \Sigma{\rm-skew\  Armendariz}\  \subsetneqq\  {\rm weak}\ \Sigma{\rm-skew\ Armendariz}
\end{equation*}
\begin{equation*}
(\Sigma, \Delta){\rm -Armendariz} \ \subsetneqq\ \Sigma{\rm-skew\  Armendariz}
\end{equation*}
\begin{equation*}
(\Sigma, \Delta){\rm -weak\ Armendariz}\ \subsetneqq\  {\rm weak}\ \Sigma{\rm-skew\ Armendariz}
\end{equation*}
As we can appreciate, the more general class of rings consists of the weak $\Sigma$-skew Armendariz. The purpose of this paper is to introduce a generalization of these rings, the {\em weak skew-Armendariz} rings (Definition \ref{Definition2.42008} below).
\section{Skew-Armendariz and weak-Armendariz rings}\label{skewArmendarizsection}
In this section we introduce two new notions of Armendariz for $\sigma$-PBW extensions: skew-Armendariz (Definition \ref{perritogatito}) and a more general notion, the weak skew-Armendariz (Definition \ref{Definition2.42008}). These definitions generalize the treatments developed for both classical polynomial rings and Ore extensions of injective type (c.f. \cite{Armendariz1974}, \cite{RegeChhawchharia1997}, \cite{AndersonCamillo1998}, \cite{HongKimKwak2000}, \cite{KimLee2000}, \cite{Hirano2002}, \cite{Huhetal2002}, \cite{HongKimKwak2003}, \cite{LeeWong2003}, \cite{Matczuk2004}, \cite{MoussaviHashemi2005},  \cite{MohamedMoufadal2006}, and \cite{NasrMoussavi2007}, and  \cite{NasrMoussavi2008}, \cite{Reyes2015}). We show also that the families of Armendariz rings defined in \cite{NinoReyes2016} and \cite{ReyesSuarez2016C} are contained in the family of skew-Armendariz and weak skew-Armendariz, but the converse  is false.
\begin{definition}\label{perritogatito}
Let $R$ be a ring and $A$ a skew PBW extension of $R$. We say that $R$ is a {\em skew-Armendariz} ring, if for polynomials $f=a_0+a_1X_1+\dotsb + a_mX_m$ and $g=b_0+b_1Y_1 + \dotsb + b_tY_t$ in $A$, $fg=0$ implies $a_0b_k=0$, for each $0\le k\le t$.
\end{definition}
Note that every Armendariz ring is skew-Armendariz, where $\sigma_i={\rm id}_R$ and $\delta_i=0$ ($1\le i\le n$), and every $\Sigma$-skew Armendariz ring is also a skew-Armendariz ring. If $R$ is $\Sigma$-rigid, the elements $c_{i,j}$ are invertible (Definition \ref{gpbwextension} (iv)), and they are at the center of $R$, from \cite{Reyes2015}, Proposition 3.6 we know that $R$ is skew-Armendariz.
\begin{definition}\label{Definition2.42008}
Let $R$ be a ring and $A$ a skew PBW extension of $R$. We say that $R$ is a {\em weak skew-Armendariz} ring, if for linear polynomials $f=a_0+a_1x_1+\dotsb + a_nx_n$, and $g=b_0+b_1x_1+\dotsb + b_nx_n$ in $A$, $fg=0$ implies $a_0b_k=0$, for every $0\le k\le n$.
\end{definition}
We can see that every skew-Armendariz ring is weak skew-Armendariz.  However, a weak Armendariz ring is not necessarily Armendariz. As an illustration of this fact in the case of Ore extensions, see \cite{LeeWong2003}, Example 3.2. Of course, every weak $\Sigma$-skew Armendariz ring is a weak skew-Armendariz ring. So, we have the relations
{\small{\begin{equation*}
\Sigma{\rm-rigid}\ \subsetneqq\ (\Sigma, \Delta){\rm -Armendariz}\ \subsetneqq\ \Sigma{\rm-skew\  Armendariz} \ \subsetneqq\ {\rm skew-Armendariz}
\end{equation*}}}
and
{\small{\begin{equation*}
\Sigma{\rm-rigid}\ \subsetneqq\ (\Sigma, \Delta){\rm -weak\ Armendariz}\ \subsetneqq\ {\rm weak}\ \Sigma{\rm-skew\ Armendariz}\ \subsetneqq\ {\rm weak\ skew-Armendariz}
\end{equation*}}}
In this way, the results presented in this paper for skew-Armendariz and weak skew-Armendariz rings generalize all results established in the previous papers \cite{Reyes2015}, \cite{NinoReyes2016}, and \cite{ReyesSuarez2016C}, for skew PBW extensions, and in particular, for Ore extensions of injective type.

We start with some key results about skew Armendariz and weak skew Armendariz rings. Lemma \ref{Lemma2.52008} extends \cite{NinoReyes2016}, Proposition 3.8, and \cite{ReyesSuarez2016C}, Lemma 3.3.
\begin{lemma}\label{Lemma2.52008}
If $R$ is a weak skew-Armendariz ring, the equality $ab=0$ implies $\sigma^{\alpha}(a)\delta^{\alpha}(b) = \delta^{\alpha}(a)b=0$, for each $\alpha\in \mathbb{N}^n$.
\begin{proof}
We only show the case $\sigma_i(a)\delta_i(b)=\delta_i(a)b=0$, for $i=1,\dotsc,n$. Since $ab=0$, then $0=\delta_i(ab)=\sigma_i(a)\delta_i(b)+\delta_i(a)b$, or equivalently, $\delta_i(a)b=-\sigma_i(a)\delta_i(b)$. Let $f, g\in A$ given by  $f = \delta_i(a) + 0x_1 + \dotsb + 0x_{i-1} + \sigma_i(a)x_i + 0x_{i+1} + \dotsb + 0x_n$, and $g =  b + bx_1 + \dotsb + bx_n$, respectively. Note that $fg=0$:
\begin{align*}
fg = &\ \delta_i(a)b + \delta_i(a)bx_1 + \dotsb + \delta_i(a)bx_n + \sigma_i(a)x_ib + \sigma_i(a)x_ibx_1 + \dotsb + \sigma_i(a)x_ibx_n\\
= &\ \delta_i(a)b + \delta_i(a)bx_1 + \dotsb + \delta_i(a)
bx_n + \sigma_i(a)[\sigma_i(b)x_i + \delta_i(b)] + \sigma_i(a)[\sigma_i(b)x_i + \delta_i(b)]x_1\\
+ &\ \dotsb + \sigma_i(a)[\sigma_i(b)x_i + \delta_i(b)]x_n\\
= &\ \delta_i(a)b + \delta_i(a)bx_1 + \dotsb + \delta_i(a)bx_n + \sigma_i(a)\sigma_i(b)x_i + \sigma_i(a)\delta_i(b) + \sigma_i(a)\sigma_i(b)x_ix_1\\
+ &\ \sigma_i(a)\delta_i(b)x_1 + \dotsb + \sigma_i(a)\sigma_i(b)x_ix_n + \sigma_i(a)\delta_i(b)x_n = 0.
\end{align*}
From Definition \ref{Definition2.42008}  we obtain $\delta_i(a)b=0$, so $\sigma_i(a)\delta_i(b)=0$. \end{proof}
\end{lemma}
In \cite{Matczuk2004} and \cite{ChenTong2005}, both authors of those papers give a positive answer to the following question formulated in \cite{HongKimKwak2003}, p. 115: Let $\sigma$ be a monomorphism (or automorphism) of a (commutative) reduced ring $B$ and $B$ be a $\sigma$-skew Armendariz. Is $B$ $\sigma$-rigid? The content of Theorem \ref{parciaaaaal} is the generalization of this answer to skew PBW extensions. We suppose that the elements $c_{i,j}$ in Definition \ref{gpbwextension} (iv) are invertible and commute with every element of $R$. These assumptions are satisfied for a lot of algebras, for example: any Ore extension $R[x;\sigma,\delta]$,  additive analogue of the Weyl algebra, multiplicative analogue of the Weyl algebra, quantum algebra $\cU'(\mathfrak{so}(3,\Bbbk))$, $3$-dimensional skew polynomial algebras, Dispin algebra $\cU(osp(1,2))$, Woronowicz algebra $\cW(\mathfrak{sl}(2,\Bbbk))$, complex algebra $V_q(\mathfrak{sl}_3(\mathbb{C}))$, $q$-Heisenberg algebra, quantum enveloping algebra of $\mathfrak{sl}(2,\Bbbk),\cU_q(\mathfrak{sl}(2,\Bbbk))$, PBW extensions, almost normalizing extensions, solvable polynomial rings, diffusion algebras, $G$-algebras with $d_{i,j}$ linear, and others. It is clear that any $\sigma$-PBW extension over a field $\Bbbk$ satisfies these assumptions.
\begin{theorem}\label{parciaaaaal}
If $A$ is a skew PBW extension of a ring $R$, then the following statements are equivalent:
\begin{enumerate}
\item [\rm (i)] $R$ is reduced and skew-Armendariz;
\item [\rm (ii)] $R$ is $\Sigma$-rigid;
\item [\rm (iii)] $A$ is reduced.
\end{enumerate}
\begin{proof}
(ii) $\Leftrightarrow$ (iii) This equivalence follows from \cite{Reyes2015}, Proposition 3.5. (ii) $\Rightarrow$ (i) From \cite{Reyes2015} we know that a $\Sigma$-rigid ring is reduced, and as we saw above, every $\Sigma$-rigid ring is also skew-Armendariz. Let us see (i) $\Rightarrow$ (ii) Suppose that $R$ is reduced, skew-Armendariz and is not $\Sigma$-rigid. Then there exists $\beta\in \mathbb{N}^n$ with $a\sigma^{\beta}(a)=0$ and $a\neq 0$. Note that $\sigma^{\beta}(a)\sigma^{\beta}(\sigma^{\beta}(a)) = \sigma^{\beta}(a\sigma^{\beta}(a)) = 0$. Using that $R$ is reduced, the equality $(\sigma^{\beta}(a)a)^{2} = \sigma^{\beta}(a)a\sigma^{\beta}(a)a=0$ implies $\sigma^{\beta}(a)a=0$. Equivalently, since $a\neq 0$, $\sigma^{\beta}$ is injective, and $R$ is reduced, then $\sigma^{\beta}(a)\neq 0$ and $(\sigma^{\beta}(a))^{2}\neq 0$. With this in mind, consider the elements $f=\sigma^{\beta}(a) + \sigma^{\beta}(a)x^{\beta}$, $g=a-\sigma^{\beta}(a)x^{\beta}$. Then
\begin{align}
fg = &\ (\sigma^{\beta}(a) + \sigma^{\beta}(a)x^{\beta})(a-\sigma^{\beta}(a)x^{\beta}) \notag \\
= &\ \sigma^{\beta}(a)a - (\sigma^{\beta}(a))^{2}x^{\beta} + \sigma^{\beta}(a)x^{\beta}a - \sigma^{\beta}(a)x^{\beta}\sigma^{\beta}(a)x^{\beta} \notag \\
= &\ - (\sigma^{\beta}(a))^{2}x^{\beta} + \sigma^{\beta}(a)[\sigma^{\beta}(a)x^{\beta} + p_{\beta, a}]- \sigma^{\beta}(a)[\sigma^{\beta}(\sigma^{\beta}(a))x^{\beta} + q_{\beta, \sigma^{\beta}(a)}]x^{\beta}\notag \\
= &\ \sigma^{\beta}(a)p_{\beta, a} - \sigma^{\beta}(a\sigma^{\beta}(a))x^{\beta}x^{\beta} - \sigma^{\beta}(a)q_{\beta, \sigma^{\beta}(a)}x^{\beta}\notag \\
= &\ \sigma^{\beta}(a)p_{\beta, a} - \sigma^{\beta}(a)q_{\beta, \sigma^{\beta}(a)}x^{\beta}, \notag
\end{align}
where $p_{\beta,a}=0$ or $\deg(p_{\beta,a})<|\beta|$, if $p_{\beta, r}\neq 0$, and $q_{\beta,\sigma^{\beta}(a)}=0$ or $\deg(q_{\beta,\sigma^{\beta}(a)})<|\beta|$, if $q_{\beta,\sigma^{\beta}(a)}\neq 0$. Since $a\sigma^{\beta}(a) =\sigma^{\beta}(a) a=0$, Remark \ref{juradpr} and Lemma \ref{Lemma2.52008} guarantee that $\sigma^{\beta}(a)p_{\beta, a} = \sigma^{\beta}(a)q_{\beta, \sigma^{\beta}(a)}x^{\beta} = 0$, so $fg=0$. By assumption, $R$ is skew-Armendariz, that is, $-(\sigma^{\beta}(a))^{2}=0$, but $-(\sigma^{\beta}(a))^{2}\neq 0$, i.e., we have obtained a contradiction. Hence, $R$ is $\Sigma$-rigid.
\end{proof}
\end{theorem}
\begin{corollary}[\cite{NinoReyes2016}, Theorem 3.9, and \cite{ReyesSuarez2016C}, Theorem 3.6]
If $A$ is a skew PBW extension of a ring $R$, then the following statements are equivalent: {\rm (i)} $R$ is reduced and $(\Sigma, \Delta)$-Armendariz {\rm (}$\Sigma$-skew Armendariz{\rm )} {\rm (ii)} $R$ is $\Sigma$-rigid {\rm (iii)} $A$ is reduced.
\end{corollary}

\begin{remark}
From \cite{NinoReyes2016} we know that $\Sigma$-rigid rings $\subsetneq (\Sigma, \Delta)$-Armendariz rings; from \cite{ReyesSuarez2016C}, we know that  $\Sigma$-rigid rings $\subsetneq \Sigma$-skew Armendariz rings, and hence $\Sigma$-rigid rings $\subsetneq$ skew-Armendariz rings. Hence, \cite{NinoReyes2016}, Theorem 3.9, \cite{ReyesSuarez2016C}, Theorem 3.6, and Theorem \ref{parciaaaaal}, show that if we assume that the ring $R$ is reduced, then for $\Sigma$-rigid rings the notions of $(\Sigma, \Delta)$-Armendariz, $\Sigma$-skew Armendariz, and skew-Armendariz, coincide. This fact shows the importance of consider skew PBW extensions over non-reduced rings with the aim of obtaining ring theoretical properties more general than the established in all these papers.
\end{remark}

The next proposition generalizes \cite{ReyesSuarez2016C}, Proposition 3.8.
\begin{proposition}\label{Proposition3.12008mio}
If $R$ is a weak skew-Armendariz ring, and $e\in R$ is an idempotent element, we have $\sigma_{i}(e)=e$ and $\delta_i(e)=0$, for every $i=1,\dotsc, n$.
\begin{proof}
Consider an idempotent element $e$ of $R$. Then $\delta_i(e)=\sigma_i(e)\delta_i(e)+\delta_i(e)e$. Let $f, g\in A$ given by $f = \delta_i(e) + 0x_1 + \dotsb + 0x_{i-1} + \sigma_i(e)x_i + 0x_{i+1} + \dotsb + 0x_n$, and $g = e-1 + (e-1)x_1 + \dotsb + (e-1)x_n$, respectively.
Recall that $\delta_i(1)=0$, for every $i$. Let us show that $fg=0$:
\begin{align*}
fg = &\ \delta_i(e)(e-1) + \biggl(\sum_{j=1}^{n} \delta_i(e)(e-1)x_j\biggr) + \sigma_i(e)x_i(e-1) + \sum_{j=1}^{n} \sigma_i(e)x_i(e-1)x_j\\
= &\ \delta_i(e)(e-1) + \biggl(\sum_{j=1}^{n} \delta_i(e)(e-1)x_j\biggr) + \sigma_i(e)[\sigma_i(e-1)x_i + \delta_i(e-1)] \\
+ &\ \sum_{j=1}^{n} \sigma_i(e)[\sigma_i(e-1)x_i + \delta_i(e-1)]x_j.
\end{align*}
Equivalently,
\newpage

\begin{align*}
fg = &\ \delta_i(e)(e-1) + \biggl(\sum_{j=1}^{n} \delta_i(e)(e-1)x_j\biggr) + \sigma_i(e) [(\sigma_i(e) - \sigma_i(1))x_i + \delta_i(e)]\\
+ &\ \sum_{j=1}^{n} \sigma_i(e) [(\sigma_i(e) - \sigma_i(1))x_i + \delta_i(e)]x_j\\
= &\ \delta_i(e)(e-1) + \biggl(\sum_{j=1}^{n} \delta_i(e)(e-1)x_j\biggr) + \sigma_i(e)[\sigma_i(e)x_i - x_i +\delta_i(e)]\\
+ &\ \sum_{j=1}^{n} \sigma_i(e)[\sigma_i(e)x_i - x_i +\delta_i(e)]x_j\\
= &\ \delta_i(e)e - \delta_i(e) + \biggl(\sum_{j=1}^{n}(\delta_i(e)e - \delta_i(e))x_j\biggr) + \sigma_i(e)x_i - \sigma_i(e)x_i + \sigma_i(e)\delta_i(e)\\
+ &\ \sum_{j=1}^{n} (\sigma_i(e)x_i - \sigma_i(e)x_i + \sigma_i(e)\delta_i(e))x_j\\
= &\ \delta_i(e)e - \delta_i(e) + \sum_{j=1}^{n} \delta_i(e)ex_j - \sum_{j=1}^{n} \delta_i(e)x_j + \sigma_i(e)\delta_i(e) + \sum_{j=1}^{n} \sigma_i(e)\delta_i(e)x_j = 0.
\end{align*}
From Definition \ref{Definition2.42008} we obtain $\delta_i(e)(e-1) = 0$, i.e., $\delta_i(e)e = \delta_i(e)$, and hence $\sigma_i(e)\delta_i(e) = 0$.

Now, consider the elements $s$ and $t$ of $A$ given by $s=  \delta_i(e) - (1-\sigma_i(e))x_i$ and $t=e + \sum_{j=1}^{n} ex_j$, respectively. Let us show that $st=0$:
\begin{align*}
st = &\ \delta_i(e)e + \biggl(\delta_i(e)e\sum_{j=1}^{n} x_j\biggr) - (1-\sigma_i(e))x_ie - \biggl((1-\sigma_i(e))x_ie\sum_{j=1}^{n} x_j\biggr)\\
= &\ \delta_i(e)e + \biggl(\delta_i(e)e\sum_{j=1}^{n}x_j\biggr) - x_ie + \sigma_i(e)x_ie - x_ie\sum_{j=1}^{n} x_j + \sigma_i(e)x_ie\sum_{j=1}^{n} x_j\\
= &\ \delta_i(e)e + \biggl(\delta_i(e)e\sum_{j=1}x_j\biggr) - (\sigma_i(e)x_i + \delta_i(e)) + \sigma_i(e)(\sigma_i(e)x_i + \delta_i(e))\\
- &\ \biggl((\sigma_i(e)x_i + \delta_i(e))\sum_{j=1}^{n} x_j\biggr) + \sigma_i(e)(\sigma_i(e)x_i + \delta_i(e))\sum_{j=1}^{n} x_j\\
= &\ \delta_i(e)e + \biggl(\delta_i(e)e\sum_{j=1}^{n}x_j\biggr) - \sigma_i(e)x_i - \delta_i(e) + \sigma_i(e)x_i + \sigma_i(e)\delta_i(e) - \sigma_i(e)x_i\sum_{j=1}^{n}x_j \\
- &\ \delta_i(e)\sum_{j=1}^{n}x_j + \sigma_i(e)x_i\sum_{j=1}^{n} x_j + \sigma_i(e)\delta_i(e)\sum_{j=1}^{n}x_j.
\end{align*}
Since $\delta_i(e)=\delta_i(e)e$ and $\sigma_i(e)\delta_i(e)=0$, then $st=0$. By Armendariz condition we know that $\delta_i(e)e=0$, which shows that $\delta_i(e)=0$.

Consider the elements $u, v \in A$ given by $u=1-e + (1-e)\sigma_i(e)x_i$ and $v=e+(e-1)\sigma_i(e)x_i$. We have the equalities
\begin{align*}
uv = &\ e+ (e-1)\sigma_i(e)x_i - e^2 - e(e-1)\sigma_i(e)x_i + (1-e)\sigma_i(e)x_ie + (1-e)\sigma_i(e)x_i(e-1)\sigma_i(e)x_i\\
= &\ e\sigma_i(e)x_i - \sigma_i(e)x_i - e\sigma_i(e)x_i + e\sigma_i(e)x_i + (1-e)\sigma_i(e)(\sigma_i(e)x_i + \delta_i(e))\\
+ &\ (1-e)\sigma_i(e)(\sigma_i(e)x_i - x_i + \delta_i(e))\sigma_i(e)x_i\\
= &\ -\sigma_i(e)x_i + e\sigma_i(e)x_i + \sigma_i(e)x_i + \sigma_i(e)\delta_i(e) - e\sigma_i(e)x_i - e\sigma_i(e)\delta_i(e) \\
+ &\ [\sigma_i(e)x_i - \sigma_i(e)x_i + \sigma_i(e)\delta_i(e) - e\sigma_i(e)x_i + e\sigma_i(e)x_i - e\sigma_i(e)\delta_i(e)]\sigma_i(e)x_i = 0.
\end{align*}
Using that $\delta_i(e)=0$, we obtain $(1-e)(e-1)\sigma_i(e)=0$, i.e., $e\sigma_i(e) = \sigma_i(e)$.

Finally, let $w=e+e(1-\sigma_i(e))x_i,\ z=1-e - e(1-\sigma_i(e))x_i$ be elements of $A$. Then
\begin{align*}
wz = &\ e-e^2-e^2(1-\sigma_i(e))x_i + e(1-\sigma_i(e))x_i - e(1-\sigma_i(e))x_ie - e(1-\sigma_i(e))x_ie(1-\sigma_i(e))x_i\\
= &\ -e(1-\sigma_i(e))(\sigma_i(e)x_i + \delta_i(e)) - e(1-\sigma_i(e))[\sigma_i(e(1-\sigma_i(e)))x_i + \delta_i(e(1-\sigma_i(e)))]x_i.
\end{align*}
Using that $\delta_i(e)=0$ and $e\sigma_i(e)=\sigma_i(e)$, we can see that $wz=0$. Hence, $e(-e(1-\sigma_i(e))) = 0$, which shows that $e\sigma_i(e)=e$, and so $\sigma_i(e)=e$.
\end{proof}
\end{proposition}
Proposition \ref{Proposition3.22008} generalizes \cite{ReyesSuarez2016C}, Proposition 3.7.
\begin{proposition}\label{Proposition3.22008}
Let $R$ be a skew-Armendariz ring. If $e^2=e\in  A$, with $e=\sum_{i=0}^{m} e_iX_i$, then $e\in R$.
\begin{proof}
Let $e=e_0+e_1X_1+\dotsb  + e_mX_m$ be an element of $A$ with $e^2=e$. Since $(e_0+e_1X_1+\dotsb e_nX_n)((1-e_0) - e_1X_1 - \dotsb - e_nX_n) = ((1-e_0) - e_1X_1 - \dotsb - e_nX_n) (e_0+e_1X_1+\dotsb e_nX_n)$, the assumption on $R$ implies $e_0(1-e_0)=(1-e_0)e_i=e_0e_i$ $(1\le i\le n)=0$. Hence $e_i=0$, for every $i$, which shows that $e=e_0=e^{2}_{0}$.
\end{proof}
\end{proposition}
Proposition \ref{Theorem3.32008} extends \cite{ReyesSuarez2016C}, Proposition 3.9.
\begin{proposition}\label{Theorem3.32008}
Every weak skew-Armendariz ring is Abelian.
\begin{proof}
Let $e^2=e, a \in R$. Consider the elements $f, g$ of $A$ given by $f=e - \sum_{i=1}^{n} ea(1-e)x_i$, and $g=1-e + \sum_{i=1}^{n} ea(1-e)x_i$. Since
\begin{align*}
fg = &\ e - e^2 + e^2a(1-e) \sum_{i=1}^{n} x_i - \sum_{i=1}^{n} ea(1-e)x_i + \sum_{i=1}^{n} ea(1-e) x_i e \\
- &\ \biggl(\sum_{i=1}^{n} ea(1-e)x_iea(1-e)\biggr) \biggl(\sum_{i=1}^{n} x_i \biggr)\\
= &\ \sum_{i=1}^{n} ea(1-e)(\sigma_i(e)x_i + \delta_i(e)) - \biggl(\sum_{i=1}^{n} ea(1-e)x_i(ea-eae)\biggr) \biggl(\sum_{i=1}^{n} x_i\biggr)\\
= &\ (ea-eae)\sum_{i=1}^{n} \sigma_i(e)x_i  \\
- &\ \biggl(\sum_{i=1}^{n} (ea-eae)(\sigma_i(ea-eae)x_i + \delta_i(ea-eae))\biggr) \biggl(\sum_{i=1}^{n}x_i\biggr)
\end{align*}
\begin{align*}
fg = &\ (ea-eae)\sum_{i=1}^{n} \sigma_i(e)x_i \\
- &\ \biggl(\sum_{i=1}^{n} (ea-eae) ((e\sigma_i(a) - e\sigma_i(a)e)x_i + \delta_i(ea) - \delta_i(eae))\biggr) \biggl(\sum_{i=1}^{n} x_i\biggr) \\
= &\ ea\sum_{i=1}^{n} \sigma_i(e)x_i - eae\sum_{i=1}^{n} \sigma_i(e)x_i \\
- &\ \biggl(\sum_{i=1}^{n} (ea-eae) (e\sigma_i(a)x_i - e\sigma_i(a)ex_i + \sigma_i(e)\delta_i(a) - \sigma_i(ea)\delta_i(e) - \delta_i(ea)e)\biggr)\biggl(\sum_{i=1}^{n} x_i\biggr) \\
= &\ ea\sum_{i=1}^{n} \sigma_i(e)x_i - eae\sum_{i=1}^{n} \sigma_i(e)x_i \\
- &\ \biggl(\sum_{i=1}^{n} (ea-eae)(e\sigma_i(a)x_i - e\sigma_i(a)ex_i + e\delta_i(a) - e\delta_i(a)e)\biggr)\biggl(\sum_{i=1}^{n} x_i\biggr)\\
= &\ ea\sum_{i=1}^{n} \sigma_i(e)x_i - eae\sum_{i=1}^{n} \sigma_i(e)x_i -  \biggl(\sum_{i=1}^{n} eae\sigma_i(a)x_i - eae\sigma_i(a)ex_i + eae\delta_i(a) - eae\delta_i(a)e  \\
- &\ eae\sigma_i(a)x_i + eae\sigma_i(a)ex_i - eae\delta_i(a) + eae\delta_i(a)e \biggr) \biggl(\sum_{i=1}^{n} x_i\biggr) = 0.
\end{align*}
Since $R$ is weak skew-Armendariz, $eea(1-e)=0$, that is, $ea=eae$.

Now, consider the elements $p, q$ of $A$ given by $p=1-e - \sum_{i=1}^{n} (1-e)aex_i$ and $q = e + \sum_{i=1}^{n} (1-e)aex_i$. Then $pq=0$. More exactly,
\begin{align*}
pq = &\ e + \sum_{i=1}^{n} (1-e)aex_i - e^{2} - e\sum_{i=1}^{n} (1-e)aex_i - \sum_{i=1}^{n} (1-e)aex_ie \\
-&\ \biggl(\sum_{i=1}^{n} (1-e)aex_i\biggr)\biggl( \sum_{i=1}^{n} (1-e)aex_i\biggr)\\
= &\ \sum_{i=1}^{n} (1-e)aex_i - \sum_{i=1}^{n} (1-e)ae (\sigma_i(e)x_i + \delta_i(e)) - \biggl(\sum_{i=1}^{n} (1-e)aex_i(1-e)ae   \biggr) \biggl(\sum_{i=1}^{n}x_i \biggr) \\
= &\ - \biggl( \sum_{i=1}^{n} (ae-eae) (\sigma_i(ae-eae)x_i + \delta_i(ae-eae))\biggr) \biggl( \sum_{i=1}^{n}x_i \biggr)\\
= &\ - \biggl(\sum_{i=1}^{n} (ae-eae)(\sigma_i(a)ex_i - e\sigma_i(a)ex_i  + \delta_i(a)e  - \delta_i(ea)e)   \biggr)\biggl(\sum_{i=1}^{n} x_i \biggr)\\
= &\ - \biggl(\sum_{i=1}^{n} ae\sigma_i(a)ex_i -  ae\sigma_i(a)ex_i + ae\delta_i(a)e - ae\sigma_i(e)\delta_i(a)e - eae\sigma_i(a)ex_i + eae\sigma_i(a)ex_i \\
- &\ eae\delta_i(a)e + eae\sigma_i(e)\delta_i(a)e \biggr) \biggl( \sum_{i=1}^{n} x_i \biggr) = 0.
\end{align*}
By the weak skew-Armendariz condition, we know that $(1-e)(1-e)ae=0$, or equivalently, $ae=eae$. Now, as we showed above, $ea=eae$, which means that $ae=ea$, i.e., $R$ is Abelian.
\end{proof}
\end{proposition}
Propositions \ref{Proposition3.12008mio}, \ref{Proposition3.22008} and \ref{Theorem3.32008} imply the following result which generalizes \cite{ReyesSuarez2016C},  Corollary 3.10.
\begin{proposition}\label{Corollary3.42008}
If $R$ is a skew-Armendariz ring $R$, then $A$ is an Abelian ring.
\end{proposition}
\subsection{Localization of skew-Armendariz rings}\label{localizationJunior}
In this section we characterize the classical right quotient rings of skew-Armendariz rings.

Let us recall the key facts about noncommutative localization. If $B$ is a ring and $S$ is a multiplicative subset of $B$ ($1\in S$, $0\notin S$, $ss'\in S$, for every $s, s'\in S$), then the left ring of fractions of $B$ exists if and only if two conditions hold: (i) given $a\in B$ and $s\in S$ with $as=0$, there exists $s'\in S$ such that $s'a=0$; (left Ore condition) given $a\in B$ and $s\in S$, there exist $s'\in S$ and $a'\in B$ with $s'a=a's$. If these conditions hold, then the left ring of fractions of $B$ with respect to $S$ is denoted by $S^{-1}B$, and its elements are classes denoted using fractions. More exactly, $\frac{a}{s}:=\frac{b}{t}$ are equal if and only if there exist $c, d\in B$ such that $ca=db,\ cs=dt\in S$; $\frac{a}{s} + \frac{b}{t}:=\frac{ca+db}{u}$, where $u:=cs=dt\in S$, for some $c,d\in B$; $\frac{a}{s}\frac{b}{t}:=\frac{cb}{us}$, where $ua=ct$, for some $u\in S$ and $c\in B$. Similarly, it is defined the right Ore condition and hence the ring of fractions of $B$. The nonzero divisors elements of $B$ are called {\em regular} and the set of regular elements of $B$ is denoted by $S_0(B)$.  Recall that if $B$ is both left and right Ore, then its classical left ring of quotients $Q_{cl}^{l}(B)$ and its classical right ring of quotients $Q_{cl}^{r}(B)$ coincide, and it is denoted by $Q(B)$. A key result about the classical ring of quotients of $B$ is the common denominator property: if $B$ is a ring, $S\subset B$ is a multiplicative subset and $S^{-1}B$ exists, then any finite set $\{q_1,\dotsc, q_n\}$ of elements of $S^{-1}B$ posses a common denominator, i.e., there exist $r_1,\dotsc,r_n\in B$ and $s\in S$ such that $q_i=\frac{r_i}{s}$ for every $i$ (see \cite{Jategaonkar1986} for a detailed treatment of localization in non-commutative rings).
\begin{proposition}[\cite{Lezamaetal2013}, Lemma 2.6]
Let $A$ be a bijective skew PBW extension of a ring $R$. If $S\subseteq S_0(R)$ is a multiplicative subset of $R$ with $\sigma_i(S)=S$, for every $i=1,\dotsc, n$, then
\begin{enumerate}
\item [\rm (i)] If $S^{-1}R$ exists, then $S^{-1}A$ exists and it is a bijective skew PBW extension of $S^{-1}R$, denoted $S^{-1}A=\sigma(S^{-1}R)\langle x_1',\dotsc, x_n'\rangle$, where $x_i':=\frac{x_i}{1}$, and the systems of constants of $S^{-1}R$ is given by $c'_{i,j}=\frac{c_{i,j}}{1},\ c'_{i, \frac{r}{s}}:=\frac{\sigma_i(r)}{\sigma_i(s)}$, for $1\le i, j\le n$. The automorphisms $\overline{\sigma_i}$ of $S^{-1}R$ and the $\overline{\sigma_i}$-derivations $\overline{\delta_i}$ $(1\le i\le n)$, are defined by $\overline{\sigma_i}(\frac{a}{s}):=\frac{\sigma_i(a)}{\sigma_i(s)}$, and $\overline{\delta_i}(\frac{a}{s}):=-\frac{\delta_i(s)}{\sigma_i(s)}\frac{a}{s} + \frac{\delta_i(a)}{\sigma_i(s)}$. Let $\overline{\Sigma}:=\{\overline{\sigma_1}, \dotsc, \overline{\sigma_n}\}$ and $\overline{\Delta}:=\{\overline{\delta_1}, \dotsc, \overline{\delta_n}\}$.
\item [\rm (ii)] If $RS^{-1}$ exists, then $AS^{-1}$ exists and it is a bijective skew PBW extension of $RS^{-1}$, denoted $AS^{-1}=\sigma(RS^{-1})\langle x_1'',\dotsc, x_n''\rangle$, where $x_i'':=\frac{x_i}{1}$, and the systems of constants of $S^{-1}R$ is given by $c''_{i,j}=\frac{c_{i,j}}{1},\ c''_{i, \frac{r}{s}}:=\frac{\sigma_i(r)}{\sigma_i(s)}$, for $1\le i, j\le n$. The automorphisms $\overline{\sigma_i}$ of $S^{-1}R$ and the $\overline{\sigma_i}$-derivations $\overline{\delta_i}$ $(1\le i\le n)$, are defined by $\overline{\sigma_i}(\frac{a}{s}):=\frac{\sigma_i(a)}{\sigma_i(s)}$, and $\overline{\delta_i}(\frac{a}{s}):=-\frac{\sigma_i(a)}{\sigma_i(s)}\frac{\delta_i(s)}{s} + \frac{\delta_i(a)}{s}$. Let $\overline{\Sigma}:=\{\overline{\sigma_1}, \dotsc, \overline{\sigma_n}\}$ and $\overline{\Delta}:=\{\overline{\delta_1}, \dotsc, \overline{\delta_n}\}$.
\end{enumerate}
\end{proposition}
If no confusion arises, we simply denote $x_i'$ and $x_i''$ by $x_i$ for $1\le i\le n$. Now, analogously to the definitions of $\Sigma$-rigid, skew-Armendariz and weak skew-Armendariz, we consider the notions of $\overline{\Sigma}$-rigid, skew-Armendariz and weak-skew Armendariz, for the classical quotient ring $Q(R)$ of $R$.

The next theorem generalizes \cite{KimLee2000}, Theorem 16, \cite{NasrMoussavi2007}, Theorem 2.3, \cite{NasrMoussavi2008}, Theorem 4.7, and \cite{ReyesSuarez2016C}, Theorem 4.2.
\begin{theorem}\label{Theorem4.72008}
Let $A$ be a bijective skew PBW extension of a ring $R$. If the classical ring of quotients  $Q(R)$ of $R$ exists, then $R$ is weak skew-Armendariz if and only if $Q(R)$ is weak skew-Armendariz.
\begin{proof}
It is clear that if $Q(R)$ is weak skew-Armendariz, then $R$ is weak skew-Armendariz.

Conversely, consider $f=c_0^{-1}a_0 + \sum_{i=1}^{n}c_i^{-1}a_ix_i$ and $g=s_0^{-1}b_0 + \sum_{j=1}^{n}s_j^{-1}b_jx_j$ elements of $S^{-1}A$ such that $fg=0$. Let us prove that $c_0^{-1}a_0s_j^{-1}b_j = 0$, for $0\le j\le n$.

We know that there exist $a_i', b_j'\in R$ and $c, s\in S_0(R)$ satisfying $c_i^{-1}a_i=c^{-1}a_i'$ and $s_j^{-1}b_j=s^{-1}b_j'$ for $0\le i, j \le n$. In this way, we can write
\begin{align*}
0 = &\ \biggl(c^{-1}a_0' + \sum_{i=1}^{n}c^{-1}a_i'x_i\biggr)\biggl(s^{-1}b_0' + \sum_{j=1}^{n}s^{-1}b_j'x_j\biggr)\\
= &\ \biggl(a_0' + \sum_{i=1}^{n} a_i'x_i\biggr)s^{-1}\biggl(b_0' + \sum_{j=1}^{n}b_j'x_j\biggr)\\
= &\ \biggl(a_0's^{-1} + \sum_{i=1}^{n}a_i'[\overline{\sigma_i}(s^{-1})x_i + \overline{\delta_i}(s^{-1})]\biggr) \biggl(b_0' + \sum_{j=1}^{n}b_j'x_j\biggr)\\
= &\ \biggl(a_0's^{-1} + \sum_{i=1}^{n}(a_i'\sigma_i(s)^{-1}x_i - a_i'\sigma_i(s)^{-1}\delta_i(s)s^{-1})\biggr)\biggl(b_0' + \sum_{j=1}^{n}b_j'x_j\biggr).
\end{align*}
There exist $d_i\in R\ (1\le i\le n)$ and $s_2\in S_0(R)$ such that $\delta_i(s)s^{-1}=s_2^{-1}d_i$, which shows that
\[
0 = \biggl(a_0's^{-1} + \sum_{i=1}^{n} (a_i'\sigma_i(s)^{-1}x_i - a_i'\sigma_i(s)^{-1}s_2^{-1}d_i)\biggr) \biggl(b_0'+\sum_{j=1}^{n} b_j'x_j\biggr).
\]
Since there exist $a_0'', a_i'', a_i'''\in R\ (1\le i\le n)$ and $s_3, s_4, s_5\in S_0(R)$ with the relations $a_0's^{-1} = s_3^{-1}a_0'',\ a_i'\sigma_i(s)^{-1}= s_4^{-1}a_i{''}$, and $a_i'\sigma_i(s)^{-1}s_2^{-1}=s_5^{-1}a_i'''$, we have
\[
0 = \biggl(s_3^{-1}a_0'' + \sum_{i=1}^{n}(s_4^{-1}a_i''x_i - s_5^{-1}a_i'''d_i)\biggr) \biggl(b_0' + \sum_{j=1}^{n} b_j'x_j\biggr).
\]
Again, there exist $d_0, d_i', d_i''\in R\ (1\le i\le n)$ and $t\in S_0(R)$ with $s_3^{-1}a_0'' = t^{-1}d_0,\ s_4^{-1}a_i'' = t^{-1}d_i'$, and $s_5^{-1}a_i''' = t^{-1}d_i''$, which guarantees that
\begin{align}
0 = &\ \biggl(t^{-1}d_0 + \sum_{i=1}^{n} (t^{-1}d_i'x_i - t^{-1}d_i''d_i)\biggr) \biggl(b_0' + \sum_{j=1}^{n}b_j'x_j\biggr) \notag   \\
= &\ t^{-1} \biggl(d_0 + \sum_{i=1}^{n} (d_i'x_i - d_i''d_i)\biggr) \biggl(b_0' + \sum_{j=1}^{n}b_j'x_j\biggr)  \notag \end{align}
\begin{align}
0 = &\ \biggl(d_0 - \sum_{i=1}^{n}d_i''d_i +
\sum_{i=1}^{n}d_i'x_i\biggr) \biggl(b_0' +
\sum_{j=1}^{n}b_j'x_j\biggr). \label{jjjjjj}
\end{align}
By the Armendariz condition on $R$, from (\ref{jjjjjj}) we obtain the relations given by
\begin{equation}\label{livernee}
\biggl(d_0 - \sum_{i=1}^{n} d_i''d_i  \biggr)b_j' = 0\ \ \ (0\le j\le n),\ {\rm and}\ \  d_i'\sigma_i(b_j') = 0\ \ (1\le i\le n)\ (0\le j\le n).
\end{equation}
From the reasoning above we have the following equivalences for expressions in (\ref{livernee}):
\begin{align}
\biggl(d_0 - \sum_{i=1}^{n} d_i''d_i\biggr)b_j' = 0 &\ \Leftrightarrow t^{-1}\biggl(d_0 - \sum_{i=1}^{n} d_i''d_i\biggr)b_j' = 0 \notag \\
&\ \Leftrightarrow \biggl(t^{-1}d_0 - \sum_{i=1}^{n} t^{-1}d_i''d_i\biggr)b_j' = 0 \notag \\
&\ \Leftrightarrow
\biggl(s_3^{-1}a_0'' - \sum_{i=1}^{n} s_5^{-1}a_i'''d_i\biggr)b_j' = 0 \notag \\
&\ \Leftrightarrow \biggl(a_0's^{-1} - \sum_{i=1}^{n}a_i'\sigma_i(s)^{-1}s_2^{-1}d_i\biggr)b_j' = 0 \notag \\
&\ \Leftrightarrow \biggl(a_0's^{-1} - \sum_{i=1}^{n} a_i'\sigma_i(s)^{-1}\delta_i(s)s^{-1}\biggr)b_j'=0 \notag \\
&\ \Leftrightarrow \biggl(a_0's^{-1} +
\sum_{i=1}^{n}a_i'\overline{\delta_i}(s^{-1})\biggr)b_j' = 0\ \ \
\ (0\le j\le n), \label{zzzz}
\end{align}
and
\begin{align}
d_i'\sigma_i(b_j') = 0 &\ \Leftrightarrow  t^{-1}d_i'\sigma_i(b_j') = 0 \Leftrightarrow s_4^{-1}a_i''\sigma_i(b_j') = 0\notag \\
&\ \Leftrightarrow a_i'\sigma_i(s)^{-1}\sigma_i(b_j') = 0 \Leftrightarrow a_i'\overline{\sigma_i}(s^{-1})\sigma_i(b_j') = 0 \notag \\
&\ \Leftrightarrow a_i'\overline{\sigma_i}(s^{-1})\overline{\sigma_i}(b_j') = 0 \ \ \ \ \ (1\le i\le n)\ (0\le j\le n).\label{sevil}
\end{align}
Expression (\ref{jjjjjj}) is equivalent to
\begin{align}
0 = &\ \biggl(d_0 - \sum_{i=1}^{n} d_i''d_i\biggr)b_0' + \biggl(d_0 - \sum_{i=1}^{n}d_i''d_i\biggr)\biggl(\sum_{j=1}^{n}b_j'x_j\biggr)\notag \\
+ &\ \biggl(\sum_{i=1}^{n}d_i'x_i\biggr)b_0' + \biggl(\sum_{i=1}^{n}d_i'x_i\biggr)\biggl(\sum_{j=1}^{n}b_j'x_j\biggr)\notag  \\
= &\ \biggl(d_0 - \sum_{i=1}^{n} d_i''d_i\biggr)b_0' + \biggl(d_0 - \sum_{i=1}^{n}d_i''d_i\biggr)\biggl(\sum_{j=1}^{n}b_j'x_j\biggr) \notag  \\
+ &\
\sum_{i=1}^{n}d_i'[\sigma_i(b_0')x_i + \delta_i(b_0')] + \biggl(\sum_{i=1}^{n}d_i'x_i\biggr)\biggl(\sum_{j=1}^{n}b_j'x_j\biggr)\notag
\end{align}
\begin{align}
0 = &\ \biggl(d_0 - \sum_{i=1}^{n} d_i''d_i\biggr)b_0' + \sum_{j=1}^{n}\biggl(d_0 - \sum_{i=1}^{n}d_i''d_i\biggr)b_j'x_j + \sum_{i=1}^{n}d_i'\sigma_i(b_0')x_i + \sum_{i=1}^{n} d_i'\delta_i(b_0') \notag  \\
+ &\ \sum_{i=1}^{n} d_i'x_ib_i'x_i + \sum_{i, j\in \{1,\dotsc,n\},\ i\neq j} d_i'x_ib_j'x_j\notag\\
= &\ \biggl(d_0 - \sum_{i=1}^{n}d_i''d_i\biggr)b_0' +\sum_{i=1}^{n}d_i'\delta_i(b_0') + \sum_{j=1}^{n} \biggl\{\biggl(d_0 - \sum_{i=1}^{n}d_i''d_i\biggr)b_j' + d_j'\sigma_j(b_0')\biggr\}x_j\notag  \\
+ &\ \sum_{i=1}^{n} d_i'[\sigma_i(b_i')x_i+\delta_i(b_i')]x_i + \sum_{i, j\in \{1,\dotsc,n\},\ i\neq j} d_i'x_ib_j'x_j\notag  \\
= &\ \biggl(d_0 - \sum_{i=1}^{n}d_i''d_i\biggr)b_0' +\sum_{i=1}^{n}d_i'\delta_i(b_0') + \sum_{j=1}^{n} \biggl\{\biggl(d_0 - \sum_{i=1}^{n}d_i''d_i\biggr)b_j' + d_j'\sigma_j(b_0')\biggr\}x_j\notag  \\
+ &\ \sum_{i=1}^{n} d_i'\sigma_i(b_i')x_i^{2} + \sum_{i=1}^{n}d_i'\delta_i(b_i')x_i +  \sum_{i, j\in \{1,\dotsc,n\},\ i\neq j} d_i'x_ib_j'x_j\notag  \\
= &\ \biggl(d_0 - \sum_{i=1}^{n}d_i''d_i\biggr)b_0' +\sum_{i=1}^{n}d_i'\delta_i(b_0') + \sum_{j=1}^{n} \biggl\{\biggl(d_0 - \sum_{i=1}^{n}d_i''d_i\biggr)b_j' + d_j'\sigma_j(b_0') + d_j'\delta_j(b_j')\biggr\}x_j\notag  \\
+ &\ \sum_{i=1}^{n} d_i'\sigma_i(b_i')x_i^{2} + \sum_{i, j\in \{1,\dotsc,n\},\ i\neq j}d_i'x_ib_j'x_j \label{fererep}.
\end{align}
With the purpose of computing the last sum in (\ref{fererep}), consider the Remark \ref{notesondefsigampbw} (iii). Then
\begin{align}
0 = &\ \biggl(d_0 - \sum_{i=1}^{n}d_i''d_i\biggr)b_0' +\sum_{i=1}^{n}d_i'\delta_i(b_0') \notag \\
+ &\ \sum_{j=1}^{n} \biggl\{\biggl(d_0 - \sum_{i=1}^{n}d_i''d_i\biggr)b_j' + d_j'\sigma_j(b_0') + d_j'\delta_j(b_j')\biggr\}x_j\notag \\
+ &\ \sum_{i=1}^{n} d_i'\sigma_i(b_i')x_i^{2}\notag \\
+ &\ \sum_{i, j\in \{1,\dotsc,n\},\ i\neq j} \biggr\{[d_i'\sigma_i(b_j') + d_j'\sigma_j(b_i')c_{i,j}]x_ix_j + [d_j'\delta_j(b_i') + d_j'\sigma_j(b_i')r_i^{(i,j)}]x_i\notag \\
+ &\ [d_i'\delta_i(b_j') + d_j'\sigma_j(b_i')]x_j
+ d_j'\sigma_j(b_i')r^{(i,j)} + d_j'\sigma_j(b_i')\sum_{k=1,\ k\neq i, j}^{n} r^{(i,j)}_kx_k\biggr\}. \label{pqpqpqpq}
\end{align}
By (\ref{livernee}), the expression (\ref{pqpqpqpq}) takes the form
\[
0 = \sum_{i=1}^{n} d_i'\delta_i(b_0') + \sum_{j=1}^{n}d_j'\delta_j(b_j')x_j + \sum_{i,j\in \{1,\dotsc,n\},\ i\neq j}[d_j'\delta_j(b_i')x_i + d_i'\delta_i(b_j')x_j]
\]
or equivalently,
\[
0 = \sum_{i=1}^{n} d_i'\delta_i(b_0')  + \sum_{i,j\in \{1,\dotsc,n\}} d_i'\delta_i(b_j')x_j.
\]
By degree relations, we have necessarily the equalities
\begin{equation}\label{stem}
 \sum_{i=1}^{n} d_i'\delta_i(b_0') = \sum_{i=1}^{n} d_i'\delta_i(b_j') = 0.
\end{equation}
Indeed, for the first sum in (\ref{stem}) we note that
\begin{align}
\sum_{i=1}^{n} d_i'\delta_i(b_0') = 0 &\ \Leftrightarrow \sum_{i=1}^{n} t^{-1}d_i'\delta_i(b_0') = 0\Leftrightarrow \sum_{i=1}^{n} s_4^{-1}a_i''\delta_i(b_0') = 0\notag \\
&\ \Leftrightarrow \sum_{i=1}^{n} a_i'\sigma_i(s)^{-1}\delta_i(b_0') = 0 \Leftrightarrow \sum_{i=1}^{n} c^{-1}a_i'\sigma_i(s)^{-1}\delta_i(b_0') = 0\notag \\
&\ \Leftrightarrow \sum_{i=1}^{n} c_i^{-1}a_i\sigma_i(s)^{-1}\delta_i(b_0') = 0.\notag
\end{align}
Similarly, for the second sum in (\ref{stem}) we have
\begin{align}
\sum_{i=1}^{n} d_i'\delta_i(b_j') = 0 &\ \Leftrightarrow \sum_{i=1}^{n} t^{-1}d_i'\delta_i(b_j') = 0 \Leftrightarrow \sum_{i=1}^{n} s_4^{-1}a_i''\delta_i(b_j') = 0 \notag \\
&\ \Leftrightarrow \sum_{i=1}^{n} a_i\sigma_i(s)^{-1}\delta_i(b_j') = 0 \Leftrightarrow  \sum_{i=1}^{n} a_i\overline{\sigma_i}(s^{-1}) \overline{\delta_i}(b_j') = 0\ \ \ (0\le j\le n). \label{gggg}
\end{align}
Now, consider the elements $h = a_0'+\sum_{i=1}^{n} a_i'x_i$, and $k=s^{-1}b_j'$ with $j=1,\dotsc,n$. Then
\begin{align}
hk = &\ a_0's^{-1}b_j' + \sum_{i=1}^{n} a_i'x_is^{-1}b_j' \notag \\
= &\ a_0's^{-1}b_j' + \sum_{i=1}^{n} a_i'[\overline{\sigma_i}(s^{-1}) \overline{\sigma_i}(b_j')x_i + \overline{\delta_i}(b_j')]\notag \\
= &\ a_0s^{-1}b_j' + \sum_{i=1}^{n} a_i'\overline{\sigma_i}(s^{-1})\overline{\sigma_i}(b_j')x_i + \sum_{i=1}^{n}a_i'\overline{\delta_i}(s^{-1}b_j')\\
= &\ a_0s^{-1}b_j' + \sum_{i=1}^{n} a_i'\overline{\sigma_i}(s^{-1})\overline{\sigma_i}(b_j')x_i + \sum_{i=1}^{n} a_i'\overline{\sigma_i}(s^{-1})\overline{\delta_i}(b_j') + \sum_{i=1}^{n} a_i'\overline{\delta_i}(s^{-1})b_j'.
\end{align}
By (\ref{zzzz}), (\ref{sevil}), and (\ref{gggg}), we obtain
\[
hk = a_0s^{-1}b_j' + \sum_{i=1}^{n} a_i'\delta_i(s^{-1})b_j' = \biggl(a_0's^{-1}+\sum_{i=1}^{n} a_i'\delta_i(s^{-1}) \biggr)b_j' = 0,\ \ \ 1\le j \le n.
\]
Since there exist $m_j\in R$ for every $j$, and $n\in S_0(R)$ such that $s^{-1}b_j'=m_jn^{-1}$, then we obtain $(a_0' + \sum_{i=1}^{n} a_i'x_i)m_jn^{-1} = 0$,  whence $(a_0' + \sum_{i=1}^{n} a_i'x_i)m_j = 0$. By the Armendariz condition on $R$, $a_0'm_j=0$, for every  $0\le j\le n$. Note that $a_0'm_j=0$ is equivalent to $a_0'm_jn^{-1
}=0$, that is, $a_0's^{-1}b_j'=0$.  Since we have the equivalences
\begin{equation}
a_0's^{-1}b_j' = 0 \Leftrightarrow c^{-1}a_0's_j^{-1}b_j = 0 \Leftrightarrow c_0^{-1}a_0s_j^{-1}b_j = 0,\ \ \ 0\le j\le n,
\end{equation}
we have proved that the ring $Q(R)$ of $R$ is weak skew-Armendariz.
\end{proof}
\end{theorem}
\section{Invariant ideals, Baer, quasi-Baer, p.p. and p.q.-Baer rings}\label{QuasiBaerrings}
Kaplansky \cite{Kaplansky1968} defined a ring $B$ as a {\em Baer} (resp. {\em quasi-Baer}, which was defined by Clark in \cite{Clark1967}) ring if the right annihilator of every nonempty subset (resp. ideal) of $B$ is generated by an idempotent (the objective of these rings is to abstract various properties of von Neumann algebras and complete $\ast$-regular rings; in \cite{Clark1967}, it was used the quasi-Baer concept to characterize when a finite-dimensional algebra with unity over an algebraically closed field is isomorphic to a twisted matrix units semigroup algebra). Another generalization of Baer rings are the p.p.-rings. A ring $B$ is called {\em right} (resp. {\em left}) {\em p.p} if the right (resp. {\em left}) annihilator of each element of $B$ is generated by an idempotent (or equivalently, rings in which each principal right (resp. {\em left}) ideal is projective).  Birkenmeir et. al., \cite{BirkenmeierKimPark2001a} defined a ring {\em right} (resp. {\em left}) {\em principally quasi-Baer} (or simply {\em right} (resp. {\em left}) {\em p.q-Baer}) ring if the right annihilator of each principal right (resp. left) ideal of $B$ is generated by an idempotent. Note that in a reduced ring $B$, $B$ is Baer (resp. p.p.-) if and only if $B$ is quasi-Baer (resp. p.q.-Baer).\\

We have studied the uniform dimension (also known as Goldie dimension) and the prime ideals of $\sigma$-PBW extensions (see  \cite{Reyes2014} and \cite{LezamaAcostaReyes2015}, respectively), using the notion of invariant ideal. With this in mind and with the purpose of establishing some relations between these concepts and the property of being quasi-Baer, in this section we consider the notion of $(\Sigma, \Delta)$-quasi-baer rings. We will say that a ring $B$ with a family of automorphisms $\Sigma$ and a family of $\Sigma$-derivations $\Delta$ is called $(\Sigma,\Delta)$-{\em quasi-Baer} if the right annihilator of every $(\Sigma,\Delta)$-{\em ideal} (i.e., an ideal $I$ such that $\sigma_i(I)=I$ and $\delta_i(I)\subseteq I$, for all $i$ with $\sigma_i\in \Sigma$, and $\delta_i\in \Delta$) of $B$ is generated by an idempotent of $B$ for every $1\le i\le n$. The results presented in this section generalize the treatment developed in \cite{NasrMoussavi2008} for Ore extensions of injective type.\\

Since we concern with the properties of being Baer (quasi-Baer) over a skew PBW extension of a Baer (quasi-Baer) ring, we need to establish a criterion which allows us to extend the family $\Sigma$ of injective endomorphisms,  and the family of $\Sigma$-derivations $\Delta$ of the ring $R$ (Proposition \ref{sigmadefinition}) to the ring $A$. With this aim, for the next result consider the injective endomorphisms $\sigma_i\in \Sigma$, and the $\sigma_i$-derivations $\delta_i\in \Delta$ $(1\le i\le n)$ formulated in Proposition \ref{sigmadefinition}.
\begin{theorem}\label{omaaaaa}
Let $A$ be a skew PBW  extension of a ring $R$. Suppose that $\sigma_i\delta_j=\delta_j\sigma_i,\ \delta_i\delta_j=\delta_j\delta_i$, and $\delta_k(c_{i,j}) = \delta_k(r_l^{(i,j)}) = 0$, for $1\le i, j, l\le n$, where $c_{i,j}$ and $r_l^{(i,j)}$ are the elements established in Definition \ref{gpbwextension}. If $\overline{\sigma_{k}}:A\to A$ and $\overline{\delta_k}:A\to A$ are the functions given by $\overline{\sigma_{k}}(f):=\sigma_k(a_0)+\sigma_k(a_1)X_1 + \dotsb + \sigma_k(a_m)X_m$ and $\overline{\delta_k}(f):=\delta_k(a_0) + \delta_k(a_1)X_1 + \dotsb + \delta_k(a_m)X_m$, for every $f=a_0 + a_1X_1+\dotsb + a_mX_m\in A$, respectively, and $\overline{\sigma_k}(r):=\sigma_i(k)$, for every $1\le i\le n$, then $\overline{\sigma_k}$ is an injective endomorphism of $A$ and $\overline{\delta_k}$ is a $\overline{\sigma_k}$-derivation of $A$.
\begin{proof}
It is clear that $\overline{\sigma_i}$ is an injective endomorphism of $A$, and that $\overline{\delta_i}$ is an additive map of $A$, for every $1\le i\le n$. Next, we show that $\overline{\delta_i}(fg)=\overline{\sigma_i}(f)\overline{\delta_i}(g) + \overline{\delta_i}(f)g$, for $f,g\in A$.

Consider the elements $f=a_0 + a_1X_1+a_2X_2+\dotsb + a_mX_m$ and $g=b_0 + b_1Y_1+b_2Y_2+\dotsb + b_tY_t$. Since $\overline{\sigma_k}$ and $\overline{\delta_k}$ are additive, for every $i$, it is enough to show that
\begin{align}\label{serpi}
\overline{\delta_k}(a_iX_ib_jY_j) = \overline{\sigma_k}(a_iX_i)\overline{\delta_k}(b_jY_j) + \overline{\delta_k}(a_iX_i)b_jY_j,
\end{align}
for every $1\le i, j \le n$. As an illustration of the necessity of the assumptions above, consider the next particular computations:
\begin{align}
\overline{\delta_k}(bx_jax_i) = &\ \overline{\delta_k}b(\sigma_j(a)x_j + \delta_j(a))x_i) = \overline{\delta_k}(b\sigma_j(a)x_jx_i + b\delta_j(a)x_1)\notag \\
= &\ \overline{\delta_k}\biggl(b\sigma_j(a)\biggl(c_{i,j}x_ix_j + r_0 + \sum_{l=1}^{n}r_lx_l\biggr) + b\delta_j(a)x_i\biggr)\notag \\
= &\ \overline{\delta_k}\biggl(b\sigma_j(a)c_{i,j}x_ix_j + b\sigma_j(a)r_0 + b\sigma_j(a)\sum_{l=1}^{n}r_lx_l + b\delta_j(a)x_i\biggr) \notag \\
= &\ \delta_k(b\sigma_j(a)c_{i,j})x_ix_j + \delta_k(b\sigma_j(a)r_0) + \sum_{l=1}^{n}\delta_k(b\sigma_j(a)r_l)x_l + \delta_k(b\delta_j(a))x_i\notag \\
= &\ \sigma_k(b\sigma_j(a))\delta_j(c_{i,j})x_ix_j + \delta_k(b\sigma_j(a))c_{i,j}x_ix_j + \sigma_k(b\sigma_j(a))\delta_i(r_0) \notag \\
+ &\ \delta_k(b\sigma_j(a))r_0
+ \sum_{l=1}^{n}\sigma_k(b\sigma_j(a))\delta_i(r_l)x_l + \sum_{l=1}^{n} \delta_k(b\sigma_j(a))r_lx_l \notag \\
+ &\ \sigma_k(b)\delta_k(\delta_j(a))x_i + \delta_k(b)\delta_j(a)x_i\notag \\
=  &\ \sigma_k(b)\sigma_k(\sigma_j(a))\delta_j(c_{i,j})x_ix_j + \sigma_k(b)\delta_k(\sigma_j(a))c_{i,j}x_ix_j + \delta_k(b)\sigma_j(a)c_{i,j}x_ix_j \notag \\
+ &\ \sigma_k(b)\sigma_k(\sigma_j(a))\delta_i(r_0) + \sigma_k(b)\delta_k(\sigma_j(a))r_0 + \delta_k(b)\sigma_j(a)r_0 \notag \\
+ &\ \sum_{l=1}^{n} \sigma_k(b)\sigma_k(\sigma_j(a))\delta_i(r_l)x_l + \sum_{l=1}^{n} \sigma_k(b)\delta_k(\sigma_j(a))r_lx_l +\sum_{l=1}^{n}\delta_k(b)\sigma_j(a)r_lx_l \notag \\
+ &\ \sigma_k(b)\delta_k(\delta_j(a))x_i + \delta_k(b)\delta_j(a)x_i.\label{colooo}
\end{align}
On the other hand,
\begin{align}
\overline{\sigma_k}(bx_j)\overline{\delta_k}(ax_i) + \overline{\delta_k}(bx_j)ax_i = &\ \sigma_k(b)x_j \delta_k(a)x_i + \delta_k(b)x_jax_i\notag \\
= &\ \sigma_k(b)(\sigma_j(\delta_k(a))x_j + \delta_j(\delta_k(a)))x_i + \delta_k(b)(\sigma_j(a)x_j + \delta_j(a))x_i\notag\\
= &\ \sigma_k(b)\sigma_j(\delta_k(a))x_jx_i + \sigma_k(b)\delta_j(\delta_k(a))x_i + \delta_k(b)\sigma_j(a)x_jx_i\notag  \\
+ &\ \delta_k(b)\delta_j(a)x_i\notag \\
= &\ \sigma_k(b)\sigma_j(\delta_k(a))\biggl(c_{i,j}x_ix_j + r_0 + \sum_{l=1}^{n}r_lx_l\biggr) + \sigma_k(b)\delta_j(\delta_k(a))x_i \notag \\
+ &\ \delta_k(b)\sigma_j(a)\biggl(c_{i,j}x_ix_j + r_0 +
\sum_{l=1}^{n} r_lx_l\biggr) + \delta_k(b)\delta_j(a)x_i\notag
\end{align}
\begin{align}
\overline{\sigma_k}(bx_j)\overline{\delta_k}(ax_i) + \overline{\delta_k}(bx_j)ax_i = &\ \sigma_k(b)\sigma_j(\delta_k(a))c_{i,j}x_ix_j + \sigma_k(b)\sigma_j(\delta_k(a))r_0 \notag \\
+ &\ \sigma_k(b)\sigma_j(\delta_k(a))\sum_{l=1}^{n} r_lx_l\notag \\
+ &\ \sigma_k(b)\delta_j(\delta_k(a))x_i + \delta_k(b)\sigma_j(a)c_{i,j}x_ix_j + \delta_k(b)\sigma_j(a)r_0 \notag \\
+ &\ \delta_k(b)\sigma_j(a)\sum_{l=1}^{n}r_lx_l + \delta_k(b)\delta_j(a)x_i.\label{coloooo}
\end{align}
If we want that the expressions (\ref{colooo}) and (\ref{coloooo}) represent the same value, that is,
\[
\overline{\delta_k}(bx_jax_i) =  \overline{\sigma_k}(bx_j)\overline{\delta_k}(ax_i) + \overline{\delta_k}(bx_j)ax_i,\ \ \ \ \ 1\le i, j, k, \le n
\]
then we have to impose that $\sigma_i\delta_j = \delta_j\sigma_i$, $\delta_i\delta_j = \delta_j\delta_i$,  $\delta_k(c_{i,j}) = \delta_k(r_l^{(i,j)}) = 0$, for $1\le i, j, l\le n$, where $c_{i,j}$ and $r_l^{(i,j)}$ are the elements established in Definition
\ref{gpbwextension}. This justifies the assumptions in our theorem.

Now, the proof of the general case, that is, the expression (\ref{serpi}), it follows from the above reasoning and Remark \ref{juradpr}. Let us see the details. Consider the following expressions:
\begin{align}
\overline{\delta_k}(a_iX_ib_jY_j) = &\ \overline{\delta_k}(a_i(\sigma^{\alpha_i}(b_j)X_i + p_{\alpha_i, b_j})Y_j)=\overline{\delta_k}(a_i\sigma^{\alpha_i}(b_j)X_iY_j + a_ip_{\alpha_i, b_j}Y_j)\notag \\
= &\ \overline{\delta_k}(a_i\sigma^{\alpha_i}(b_j)(c_{\alpha_i, \beta_j}x^{\alpha_i+\beta_j} + p_{\alpha_i,\beta_j}) + a_ip_{\alpha_i,b_j}Y_j)\notag \\
= &\ \overline{\delta_k}(a_i\sigma^{\alpha_i}(b_j)c_{\alpha_i,\beta_j}x^{\alpha_i+\beta_j} + a_i\sigma^{\alpha_i}(b_j)p_{\alpha_i, \beta_j} + a_ip_{\alpha_i,b_j}Y_j)\notag \\
= &\ \overline{\delta_k}(a_i\sigma^{\alpha_i}(b_j)c_{\alpha_i,\beta_j})x^{\alpha_i+\beta_j} + \overline{\delta_k}(a_i\sigma^{\alpha_i}(b_j)p_{\alpha_i,\beta_j}) + \overline{\delta_k}(a_ip_{\alpha_i,b_j}Y_j)\notag\\
= &\ \sigma_k(a_i\sigma^{\alpha_i}(b_j))\delta_k(c_{\alpha_i,\beta_j})x^{\alpha_i+\beta_j} + \delta_k(a_i\sigma^{\alpha_i}(b_j))c_{\alpha_i,\beta_j}x^{\alpha_i+\beta_j}\notag \\
+ &\ \overline{\delta_k}(a_i\sigma^{\alpha_i}(b_j)p_{\alpha_i,\beta_j}) + \overline{\delta_k}(a_ip_{\alpha_i,b_j}Y_j)\notag \\
= &\ \sigma_k(a_i)\sigma_k(\sigma^{\alpha_i}(b_j))\delta_k(c_{\alpha_i,\beta_j})x^{\alpha_i+\beta_j} + \sigma_k(a_i)\delta_k(\sigma^{\alpha_i}(b_j))c_{\alpha_i,\beta_j}x^{\alpha_i+\beta_j}\notag \\
+ &\ \delta_k(a_i)\sigma^{\alpha_i}(b_j)c_{\alpha_i,\beta_j}x^{\alpha_i+\beta_j} + \overline{\delta_k}(a_i\sigma^{\alpha_i}(b_j)p_{\alpha_i,\beta_j}) + \overline{\delta_k}(a_ip_{\alpha_i,b_j}Y_j), \notag
\end{align}
and
\begin{align}
\overline{\sigma_k}(a_iX_i)\overline{\delta_k}(b_jY_j) + \overline{\delta_k}(a_iX_i)b_jY_j = &\ \sigma_k(a_i)X_i\delta_k(b_j)Y_j + \delta_k(a_i)X_ib_jY_j\notag \\
= &\ \sigma_k(a_i)(\sigma^{\alpha_i}(\delta_k(b_j))X_i + p_{\alpha_i,\delta_k(b_j)})Y_j \notag \\
+ &\ \delta_k(a_i)(\sigma^{\alpha_i}(b_j)X_i + p_{\alpha_i, b_j})Y_j\notag \\
= &\ \sigma_k(a_i)(\sigma^{\alpha_i}(\delta_k(b_j))X_iY_j) + \sigma_k(a_i)p_{\alpha_i, \delta_k(b_j)}Y_j\notag \\
+ &\ \delta_k(a_i)\sigma^{\alpha_i}(b_j)X_iY_j + \delta_k(a_i)p_{\alpha_i,b_j}Y_j\notag \\
= &\ \sigma_k(a_i)\sigma^{\alpha_i}(\delta_k(b_j))(c_{\alpha_i,\beta_j}x^{\alpha_i+\beta_j} + p_{\alpha_i,\beta_j}) + \sigma_k(a_i)p_{\alpha_i,\delta_k(b_j)}Y_j\notag \\
+ &\ \delta_k(a_i)\sigma^{\alpha_i}(b_j)(c_{\alpha_i,\beta_j}x^{\alpha_i+\beta_j} + p_{\alpha_i,\beta_j}) + \delta_k(a_i)p_{\alpha_i,b_j}Y_j\notag \\
= &\ \sigma_k(a_i)\sigma^{\alpha_i}(\delta_k(b_j))c_{\alpha_i,\beta_j}x^{\alpha_i+\beta_j} + \sigma_k(a_i)\sigma^{\alpha_i}(\delta_k(b_j))p_{\alpha_i, \beta_j} \notag \\
+ &\ \sigma_k(a_i)p_{\alpha_i, \delta_k(b_j)}Y_j + \delta_k(a_i)\sigma^{\alpha_i}(b_j)c_{\alpha_i, \beta_j}x^{\alpha_i+\beta_j}\notag \\
+ &\ \delta_k(a_i)\sigma^{\alpha_i}(b_j)p_{\alpha_i, \beta_j} +
\delta_k(a_i)p_{\alpha_i, b_j}Y_j.\notag
\end{align}
By assumption, we have the equalities $\sigma_k(a_i)\delta_k(\sigma^{\alpha_i}(b_j)) = \sigma_k(a_i)\sigma^{\alpha_i}(\delta_k(b_j))$ and $\delta_k(c_{\alpha_i,\beta_j}) = 0$, which means that we need to prove the relation
\begin{align}
\overline{\delta_k}(a_i\sigma^{\alpha_i}(b_j)p_{\alpha_i,\beta_j}) + \overline{\delta_k}(a_ip_{\alpha_i,b_j}Y_j) = &\
 \sigma_k(a_i)\sigma^{\alpha_i}(\delta_k(b_j))p_{\alpha_i, \beta_j} + \sigma_k(a_i)p_{\alpha_i, \delta_k(b_j)}Y_j \notag \\
+ &\ \delta_k(a_i)\sigma^{\alpha_i}(b_j)p_{\alpha_i, \beta_j} + \delta_k(a_i)p_{\alpha_i, b_j}Y_j.\label{metal}
\end{align}
However, note that this equality is a consequence of the linearity of $\delta_k$, Remark \ref{juradpr}, and the assumptions established in the formulation of the theorem. More precisely, using these facts we have
\begin{align}
\overline{\delta}_k(a_i\sigma^{\alpha_i}(b_j)p_{\alpha_i, \beta_j}) = &\ \overline{\sigma_k}(a_i\sigma^{\alpha_i}(b_j))\overline{\delta_k}(p_{\alpha_i,\beta_j}) + \overline{\delta_k}(a_i\sigma^{\alpha_i}(b_j))p_{\alpha_i, \beta_j}\notag \\
= &\ \sigma_k(a_i)\sigma_k(\sigma^{\alpha_i}(b_j))\overline{\delta_k}(p_{\alpha_i,\beta_j}) + \sigma_k(a_i) \delta_k(\sigma^{\alpha_i}(b_j))p_{\alpha_i, \beta_j} + \delta_k(a_i)\sigma^{\alpha_i}(b_j)p_{\alpha_i,\beta_j}\notag \\
= &\ \sigma_k(a_i) \delta_k(\sigma^{\alpha_i}(b_j))p_{\alpha_i, \beta_j} + \delta_k(a_i)\sigma^{\alpha_i}(b_j)p_{\alpha_i,\beta_j}\notag \\
= &\ \sigma_k(a_i)\sigma^{\alpha_i}(\delta_k(b_j))p_{\alpha_i,\beta_j} + \delta_k(a_i)\sigma^{\alpha_i}(b_j)p_{\alpha_i, \beta_j},\label{acdc}
\end{align}
and,
\begin{align}
\overline{\delta_k}(a_ip_{\alpha_i,b_j}Y_j) = &\ \overline{\sigma_k}(a_i)\overline{\delta_k}(p_{\alpha_i,b_j}Y_j) + \overline{\delta_k}(a_i)p_{\alpha_i,b_j}Y_j\notag \\
= &\ \sigma_k(a_i)p_{\alpha_i, \delta_k(b_j)}Y_j + \delta_k(a_i)p_{\alpha_i,b_j}Y_j,\label{marliiii}
\end{align}
where we can see that expression (\ref{metal}) is precisely the sum of (\ref{acdc}) and (\ref{marliiii}). Therefore $\overline{\delta_i}$ is a $\overline{\sigma_i}$-derivation of $A$.
\end{proof}
\end{theorem}
Birkenmeier, Kim, and Park in \cite{Birkenmeieretal2001b}, defined an idempotent $e\in B$ as left (resp. right) {\em semicentral} in $B$, if $exe=xe$ (resp. $exe=ex$), for all $x\in B$. Equivalently, $e^{2}=e\in B$ is left (resp. right) semicentral if $eB$ (resp. $Be$) is an ideal of $B$. Since the right annihilator of a right ideal is an ideal, we see that the right annihilator of a right ideal is generated by a left semicentral in a quasi-Baer ring. $\cS_l(B)$ and $\cS_r(B)$ denote the sets of all left and right semicentral idempotents of $B$, respectively. Note that $\cS_l(B) \cap \cS_r(B)=\cB(B)$, where $\cB(B)$ is the set of all central idempotents of $B$. If $B$ is a semiprime ring then $\cS_l(B)=\cS_r(B)\cB(B)$.

The next theorem generalizes \cite{NasrMoussavi2008}, Theorem 3.8.
\begin{theorem}\label{Theorem3.82008mio}
Let $A$ be a bijective skew PBW  extension of a ring $R$ with $\sigma_i(e)=e$ and $\delta_i(e)=0$, for any left semicentral idempotent $e\in R$, and every $1\le i\le n$. Suppose the conditions established in Theorem \ref{omaaaaa} hold. If $R$ is a $(\Sigma,\Delta)$-quasi Baer ring, then $A$ is a $\overline{\Sigma}$-quasi Baer ring.
\begin{proof}
Let $I$ be a $\overline{\Sigma}$-ideal of $A$. Consider $I_0$ as the set of all leading coefficients of all elements of $I$. Then $I_0$ is a left ideal of $R$. Now, since $\overline{\sigma_i}(I)=I$, then $\sigma_i(I_0)=I_0$, for every $\sigma_i\in \Sigma$. It is easy to see that $\delta_i(I_0)\subseteq I_0$, for every $i$. In this way, $I_0$ is an $(\Sigma, \Delta)$-ideal of $R$, so $r_R(I_0)=eR$ for some left semicentral idempotent $e$ of $R$. The idea is to prove that $r_A(I)=eA$. Let $f\in I$ expressed as $f=a_0 + a_1X_1+\dotsb + a_mX_m$, with ${\rm deg}(f) = |{\rm exp}(X_m)|$ and ${\rm lc}(f)=a_m$. We have $fe = (a_0 + a_1X_1+\dotsb + a_{m-1}X_{m-1})e + a_mX_me$. As we can see from Remark \ref{juradpr}, the product $X_me$ involves the elements $\sigma_j(e)$ and $\delta_k(e)$ ($1\le j, k \le n$), which are zero by assumption, so $X_me=0$. Since $fe\in I$, then ${\rm lc}(fe)\in I_0$. By the assumption on $\sigma_i(e)$ for every $1\le i\le n$, we can see that ${\rm lc}(fe)= a_{m-1}e$, i.e., $a_{m-1}e\in I_0$, so $a_{m-1}e=0$. Following this reasoning, we obtain $fe=0$,  which guarantees the inclusion $eA\subseteq r_A(I)$.

Let $g\in r_A(I),\ f\in I$. Since $I$ is a $\overline{\Sigma}$-ideal, the equality $fg=0$ implies $\sigma^{\alpha_m}(f)g = 0$, where $f=a_0 + a_1X_1 + \dotsb + a_mX_m, g= b_0+b_1Y_1 + \dotsb + b_sY_s$, and  $\alpha_m = {\rm deg}(X_m)$. Since $A$ is bijective, ${\rm lc}(\sigma^{\alpha_m}(f)g) = \sigma^{\alpha_m}(a_m)\sigma^{\alpha_m}(b_s) = \sigma^{\alpha_m}(a_mb_s)=0$, from which we conclude that $a_mb_s=0$, that is, $b_s\in r_R(I_0)=eR$, so $b_s=er'$, for some element $r'\in R$. Using that $Re=eRe$, we have $eb_s=e^{2}r' = er'=b_s$. From the equality $0=fg = f(b_0 + b_1Y_1 + \dotsb + b_{s-1}Y_{s-1}) + fb_sY_s$ and the fact $fe=0$, we conclude that $feb_sY_s=fb_sY_s=0$.

Now, since  $\sigma^{\alpha_m}(f)(b_0 + b_1Y_1 + \dotsb + b_{s-1}Y_{s-1})=0$ and ${\rm lc}(\sigma^{\alpha_m}(f)g) = \sigma^{\alpha_m}(a_m)\sigma^{\alpha_m}(b_{s-1})=0$, then $a_mb_{s-1}=0$. It is easy to see that  $b_{s-1}=eb_{s-1}$, so continuing in this way we can show  that $b_jeb_j$ for $0\le j\le s$, which proves that $r_A(I)\subseteq eA$. Therefore $r_A(I)=eA$.
\end{proof}
\end{theorem}
The next theorem generalizes \cite{NasrMoussavi2008}, Theorem 3.9.
\begin{theorem}\label{Theorem3.92008}
Let $A$ be a bijective skew PBW extension of $R$. Suppose the conditions established in Theorem \ref{omaaaaa} hold. If $R$ is skew-Armendariz, then the following statements are equivalent:
\begin{enumerate}
\item $R$ is a $(\Sigma,\Delta)$-quasi Baer ring;
\item $A$ is a $\overline{\Sigma}$-quasi Baer ring;
\item $A$ is a $(\overline{\Sigma},\overline{\Delta})$-quasi Baer ring for every extended $\overline{\alpha}_i$-derivation $\overline{\delta}_i$ of $A$.
\end{enumerate}
\begin{proof}
(1) $\Rightarrow$ (2). The assertion follows from Proposition \ref{Proposition3.12008mio} and Theorem \ref{Theorem3.82008mio}.  (2) $\Rightarrow$ (3). It is easy. (3) $\Rightarrow$ (1). Let $I$ be an $(\Sigma, \Delta)$-ideal of $R$. Then $IA$ is an $(\overline{\Sigma}, \overline{\Delta})$-ideal of $A$, so $r_A(IA) = eA$ for some idempotent $e\in A$, and by Proposition \ref{Proposition3.22008}, $e\in R$. Next we show that $r_R(I) = r_A(IA)\cap R$.

Consider an element $h$ of $IA$. Then $h=\sum_{j=1}^{s} r_jf_j$, $r_j\in I,\ f_j\in A$. Since $r\in r_R(I)$, Remark \ref{juradpr} and Lemma \ref{Lemma2.52008} guarantee that $r_jf_jr=0$, for every $j$, and hence $hr=0$. Sinc $r_A(IA)\cap R \subseteq r_R(I)$, then $r_R(I)=r_A(IA)\cap R$, and using that $r_A(IA)\cap R=eA\cap R=eR$, we conclude the desired equality.
\end{proof}
\end{theorem}

\begin{remark}It is important to say that the class of $\Sigma$-quasi Baer rings (respectively, $\Delta$-quasi Baer rings) strictly contains the class of quasi-Baer rings. This assertion can be appreciated in the case of Ore extensions of bijective type. More exactly, in \cite{Hirano2001}, Example 2, it was presented an example of an $\alpha$-quasi Baer ring $B$ which is not quasi-Baer but $B[x;\alpha]$ is a quasi-Baer ring; in \cite{Hanetal2000}, Example 2, it was provided an example of a $\delta$-quasi Baer ring which is not quasi-Baer but $B[x;\delta]$ is a quasi-Baer ring. In this way, our results generalize the treatment presented in \cite{Reyes2015} about  skew PBW extensions of quasi-Baer rings. Now, the bijectivity of the injective endomorphisms $A$ in Theorem \ref{Theorem3.92008} can not be eliminated, as we can appreciate in the particular case of Ore extensions (\cite{HongKimLee2009}, Example 2).
\end{remark}
Next, we consider the relationship between the properties of being Baer and p.p. of a ring $R$ and a skew PBW extension $A$ of $R$. Recall that a ring is Abelian and Baer if and only if  is reduced and quasi-Baer (\cite{Birkenmeieretal2013}).

The next theorem generalizes \cite{NasrMoussavi2008}, Theorem 3.11, \cite{Reyes2015}, Theorem 3.9, \cite{NinoReyes2016}, Theorem 4.1, and \cite{ReyesSuarez2016C}, Theorem 5.1.
\begin{theorem}\label{Theorem3.112008}
Let $A$ be a skew $PBW$ extension of a ring $R$. If $R$ is skew-Armendariz, then $R$ is a Baer ring if and only if $A$ is a Baer ring.
\begin{proof}
Let $R$ be a skew-Armendariz and Baer ring. From Proposition \ref{Theorem3.32008} we know that $R$ is an Abelian Baer ring, and hence $R$ is a reduced Baer ring. By Theorem  \ref{parciaaaaal}, $R$ is a $\Sigma$-rigid ring. From \cite{Reyes2015}, Theorem 3.9, it follows that $A$ is a Baer ring.

Suppose that $A$ is a Baer ring. Since $R$ is a skew-Armendariz ring, $A$ is an Abelian ring by Corollary \ref{Corollary3.42008}, so $A$ is reduced. In particular, $R$ is reduced, and hence $R$ is $\Sigma$-rigid. Therefore, $R$ is Baer (\cite{Reyes2015}, Theorem 3.9).
\end{proof}
\end{theorem}
The next theorem generalizes \cite{NasrMoussavi2008}, Theorem 3.12, \cite{Reyes2015}, Theorem 3.12,  \cite{NinoReyes2016}, Theorem 4.2, and \cite{ReyesSuarez2016C}, Theorem 5.3
\begin{theorem}\label{Theorem3.122008}
Let $A$ be a bijective skew $PBW$ extension of a ring $R$. If $R$ is skew-Armendariz, then $R$ is a p.p.-ring if and only if $A$ is a p.p.-ring
\begin{proof}
Let $R$ be a skew-Armendariz and Baer ring. We know that $R$ is an  Abelian Baer ring, so reduced. Theorem \ref{parciaaaaal}  guarantees that $R$ is $\Sigma$-rigid ring, and from \cite{Reyes2015}, Theorem 3.12, it follows that $A$ is a p.p.-ring.

Let $A$ be a p.p.- ring. By assumption, $R$ is a skew-Armendariz ring, so $A$ is an Abelian ring and hence reduced. Since $R$ is $\Sigma$-rigid, $R$ is a Baer ring (\cite{Reyes2015}, Theorem 3.12).
\end{proof}
\end{theorem}
Bell in \cite{Bell1970} defined the following: a ring $B$ is said to satisfy the $IFP$ ({\em insertion of factors property}) if $r_B(a)$ is an ideal for all $a\in B$ (sometimes, a ring with $IFP$ is also called a {\em semicommutative ring}). Reduced rings have $IFP$.  In Huh et al. \cite{Huhetal2002}, Example 2, it is presented a ring with $IFP$ over which the polynomial ring need not satisfy $IFP$; in \cite{RegeChhawchharia1997}, Proposition 4.6, one find that if $B$ is Armendariz wih $IFP$, then $B[x]$ is Armendariz with $IFP$; in  \cite{Huhetal2002} one find a ring which is Armendariz  but do not satisfy $IFP$.\\

The following theorem generalizes \cite{NasrMoussavi2008}, Theorem 3.13, and  \cite{Reyes2015}, Theorems 3.10 and 3.13, and include a partial generalization of \cite{NinoReyes2016}, Theorems 4.3, 4.4, and \cite{ReyesSuarez2016C}, Corollary 5.5.
\begin{theorem}\label{Theorem3.132008}
Let $A$ be a skew $PBW$ extension of a ring $R$. If $R$ a skew-Armendariz ring, then $R$ is a quasi Baer {\rm(}respectively, p.q.-Baer{\rm )} ring with IFP if and only if $A$ is a quasi Baer {\rm (}respectively, p.q.-Baer{\rm )} ring with IFP.
\begin{proof}
If $R$ is a quasi-Baer ring with $IFP$, then $R$ is reduced, and hence Baer  (\cite{BirkenmeierKimPark2001a}, Proposition 1.14; \cite{Birkenmeieretal2001b}, Lemma 1.9). By Theorem \ref{Theorem3.112008}, $A$ is Baer. By Corollary \ref{Corollary3.42008}, $A$ is an Abelian Baer ring and hence it is reduced. Thus $A$ has $IFP$. The converse is similar.
\end{proof}
\end{theorem}
\begin{proposition}
If $A$ is a skew $PBW$ extension of an Abelian ring $R$, then the following assertions are equivalent:
\begin{enumerate}
\item [\rm (i)] $R$ is {\rm (}weak{\rm )} skew-Armendariz;
\item [\rm (ii)] For every idempotent $e\in R$ with $\sigma_i(e)=e$ and $\delta_i(e)=0$, $eR$ and $(1-e)R$ are {\rm (}weak{\rm )} skew-Armendariz;
\item [\rm (iii)] For some idempotent $e\in R$ such that $\sigma_i(e)=e$ and $\delta_i(e)=0$, $eR$ and $(1-e)R$ are {\rm (}weak{\rm )} skew-Armendariz.
\end{enumerate}
\begin{proof}
(i) $\Rightarrow$ (ii) It is easy to see that if $R$ is a (weak) skew-Armendariz ring, then for each idempotent $e\in R,\ eR$ and $(1-e)R$ are also (weak) skew-Armendariz rings. So, the assertion follows from Proposition \ref{Proposition3.12008mio}. (ii) $\Rightarrow$ (iii) It is clear. (iii) $\Rightarrow$ (i) Let $f, g \in A$ given by $f=a_0+a_1X_1+\dotsb + a_mX_m$, $g=b_0 + b_1Y_1+\dotsb + b_tY_t$, with
with $fg=0$. Since $f=(1-e)f + ef,\ g=(1-e)g+eg$, $\sigma_i(e)=e,\ \delta_i(e)=0$, and $R$ is Abelian, then $efg=0$, so ${\rm lc}(efg)=ea_0b_j=0$ for $0\le j\le t$. In this way, $(1-e)fg=0$, that is, ${\rm lc}((1-e)fg)=(1-e)a_0b_j=0$, and so $ea_0b_j$ from which $a_0b_j=0$, for every $0\le j\le m$, i.e., $R$ is skew-Armendariz.
\end{proof}
\end{proposition}
\vspace{0.5cm}
\noindent {\bf \large Acknowledgements}

\vspace{0.4cm}
The first author is supported by Grant HERMES Code 30366, Departamento de Ma\-te\-m\'a\-ti\-cas, Universidad Nacional de Colombia, Bogot\'a, Colombia.


\end{document}